\documentclass[11pt]{amsart}
\usepackage{ amsmath, amsthm, amsfonts, hyperref, graphicx, ifpdf}
\usepackage{mathrsfs,eucal,epsfig, verbatim}
\usepackage[dvipsnames,usenames]{color}
\hypersetup{
   unicode=false,          
   pdftoolbar=true,        
   pdfmenubar=true,        
   pdffitwindow=false,     
   pdfstartview={},    
   pdftitle={},    
   pdfauthor={},     
   pdfsubject={},   
   pdfcreator={},   
   pdfproducer={}, 
   pdfkeywords={}, 
   pdfnewwindow=true,      
   colorlinks=true,       
   linkcolor=blue,          
   citecolor=blue,        
   filecolor=blue,      
   urlcolor=blue          
}

\newtheorem{theorem}{Theorem}[section]
\newtheorem{cor}[theorem]{Corollary}
\newtheorem{conj}[theorem]{Question}
\newtheorem{lem}[theorem]{Lemma}
\newtheorem{pro}[theorem]{Proposition}
\newtheorem{remark}[theorem]{Remark}
\newtheorem{Def}[theorem]{Definition}
\newtheorem*{theorem1}{Theorem \ref{main1}}
\newtheorem*{theorem2}{Theorem \ref{main2}}

\theoremstyle{definition}

\DeclareMathOperator{\AH}{AH}%

\newcommand{\Mcal}{\mathcal{M}}
\newcommand{\Mscr}{\mathscr{M}}

\newcommand{\ads}{algebraic deformation space}

\newcommand{\cl}{complex length}

\newcommand{\cmc}{constant mean curvature}

\newcommand{\eus}{Euclidean space}

\newcommand{\fhtm}{fibered hyperbolic three-manifold}
\newcommand{\foc}{fibering over the circle}

\newcommand{\hg}{hyperbolic geometry}
\newcommand{\hym}{hyperbolic metric}

\newcommand{\htm}{hyperbolic three-manifold}

\newcommand{\is}{incompressible surface}
\newcommand{\ksg}{Kleinian surface group}
\newcommand{\la}{least area}

\newcommand{\ma}{minimal annulus}
\newcommand{\maxp}{maximum principle}
\newcommand{\mf}{minimal fibration}
\newcommand{\mt}{Margulis tube}
\newcommand{\mst}{maximal solid tube}
\newcommand{\ms}{minimal surface}

\newcommand{\qf}{quasi-Fuchsian}
\newcommand{\qfg}{quasi-Fuchsian group}
\newcommand{\qfm}{quasi-Fuchsian manifold}

\newcommand{\scg}{simple closed geodesic}

\newcommand{\TS}{Teichm\"{u}ller space}
\newcommand{\Tt}{Teichm\"{u}ller theory}
\newcommand{\tm}{three-manifold}
\newcommand{\tg}{totally geodesic}
\newcommand{\uhs}{upper-half space}
\newcommand{\wrt}{with respect to}

\newcommand{\be}{\begin{equation}}
\newcommand{\ene}{\end{equation}}
\newcommand{\br}{\begin{remark}}
\newcommand{\er}{\end{remark}}
\newcommand{\bl}{\begin{lem}}
\newcommand{\el}{\end{lem}}
\newcommand{\bcor}{\begin{cor}}
\newcommand{\ecor}{\end{cor}}
\newcommand{\bpro}{\begin{pro}}
\newcommand{\epro}{\end{pro}}
\newcommand{\ben}{\begin{enumerate}}
\newcommand{\een}{\end{enumerate}}
\newcommand{\bp}{\begin{proof}}
\newcommand{\ep}{\end{proof}}
\newcommand{\bpo}{\begin{pro}}
\newcommand{\epo}{\end{pro}}
\newcommand{\beq}{\begin{equation*}}
\newcommand{\eeq}{\end{equation*}}
\newcommand{\bear}{\begin{eqnarray}}
\newcommand{\eear}{\end{eqnarray}}
\newcommand{\beqar}{\begin{eqnarray*}}
\newcommand{\eeqar}{\end{eqnarray*}}
\newcommand{\bt}{\begin{theorem}}
\newcommand{\et}{\end{theorem}}

\newcommand{\C}{\mathbb{C}}
\newcommand{\R}{\mathbb{R}}

\newcommand{\T}{\mathbb{T}(\gamma)}
\renewcommand{\H}{\mathbb{H}}
\newcommand{\Ot}{\mathbb{W}}
\newcommand{\Hscr}{\mathscr{H}}
\newcommand{\Nscr}{\mathscr{N}}
\newcommand{\Ascr}{\mathscr{A}}
\newcommand{\Bscr}{\mathscr{B}}

\newcommand{\Fcal}{\mathcal{F}}
\newcommand{\Lcal}{\mathcal{L}}

\newcommand{\Tcal}{\mathcal{T}}
\newcommand{\QF}{\mathcal{QF}}

\DeclareMathOperator{\dist}{dist}%
\DeclareMathOperator{\PSL}{PSL(2,\C)}%
\DeclareMathOperator{\Area}{Area}%
\DeclareMathOperator{\Mod}{Mod(S)}%
\DeclareMathOperator{\Length}{Length}%
\DeclareMathOperator{\Otal}{Otal}

\numberwithin{equation}{section}

\allowdisplaybreaks

\def\XXint#1#2#3{{\setbox0=\hbox{$#1{#2#3}{\int}$}
    \vcenter{\hbox{$#2#3$}}\kern-.5\wd0}}

\makeatletter
\def\@citestyle{\m@th\upshape\mdseries}
\def\citeform#1{{\bfseries#1}}
\def\@cite#1#2{{%
  \@citestyle[\citeform{#1}\if@tempswa, #2\fi]}}
\@ifundefined{cite }{%
  \expandafter\let\csname cite \endcsname\cite
  \edef\cite{\@nx\protect\@xp\@nx\csname cite \endcsname}%
}{}
\makeatother

\begin{document}

\title[Minimal Surfaces in Hyperbolic three-manifolds]
{Complex length of short curves and Minimal Fibrations of
hyperbolic $3$-Manifolds fibering over the circle}


\author{Zheng Huang}
\address[Z. ~H.]{Department of Mathematics,
The City University of New York, Staten Island, NY 10314, USA}
\address{The Graduate Center, The City University of New York,
365 Fifth Ave., New York, NY 10016, USA}
\email{zheng.huang@csi.cuny.edu}

\author{Biao Wang}
\address[B. ~W.]{Department of Mathematics and Computer Science\\
         Queensborough Community College, The City University
         of New York\\
         222-05 56th Avenue Bayside, NY 11364\\}
\email{biwang@qcc.cuny.edu}

\date{July 31, 2017}



\begin{abstract}
  We investigate the {\mst}s around short {\scg}s in {\htm}s and how the {\cl}
  of curves relates to closed   {\la} incompressible {\ms}s.
  As applications, we prove the existence of closed hyperbolic three-manifolds
  fibering over the circle which are not foliated by closed imcompressible
  minimal surfaces isotopic to the fiber.
  We also show the existence of {\qfm}s containing arbitrarily many embedded closed
  incompressible {\ms}s. Our strategy is to prove main theorems
  under natural geometric conditions on the {\cl} of closed curves on a {\fhtm},
  then by computer programs, we find explicit examples where these conditions
  are satisfied.
\end{abstract}
\maketitle

\section{Introduction}


\subsection{Motivating Questions}
As fundamental objects in differential geometry, minimal hypersurfaces in {\eus} and other
Riemannian manifolds have been extensively investigated ever since the ``Plateau Problem"
in the 1930s. We are particularly interested in the $3$-dimensional case and this paper is part
of a larger goal to understand closed incompressible {\ms}s in several different classes of
{\htm}s, their connections to {\Tt}, and the ``moduli spaces" of these {\ms}s (see \cite{GHW10, HL12, HW17}).

Throughout the paper, we denote by $S$ an oriented closed
surface of genus $g \ge 2$, and we denote $\Mcal_\psi$
or $\Mcal$ a \emph{mapping torus or {\fhtm}} with monodromy $\psi$,
which is an oriented closed {\htm} that fibers
over the circle \emph{with fiber $S$} if $\psi$ is
pseudo-Anosov. We set up the following additional notations
for the paper:

\ben
\item $\Mscr$: a \emph{\qf} manifold which is diffeomorphic to
      $S \times \R$;
\item $\Tcal(S)$: {\TS} of the surface $S$;
\item ${\QF}(S)$: the {\qf} space of $S$;
\item $\AH(S)$: the {\ads} of {\ksg}s associated to the surface $S$;
\item $\Lcal = \ell + \sqrt{-1}\theta$: the {\cl} of a {\scg}
      $\gamma$ in the {\htm}, where $\ell$ is the real length,
      and $\theta$ is the twisting angle. We always assume that
      $\ell > 0$ and $\theta \in [-\pi,\pi)$.
\item $\T$ is the {\mst} around a {\scg} $\gamma$ in a {\htm},
      whose radius is denoted by $r_0$ (see Definition ~\ref{mst}).
\een

We will study {\mst}s in metrically complete {\htm}s (without parabolics).
These tubes play fundamental roles
in the quest of determining complete (or closed) {\htm}s of small volume (see for instance \cite{Mey87, Ago02, ACS06}).
Understanding how closed incompressible {\la} {\ms}s interact with deep tubes enables us to pursue
some natural questions in {\hg}. Our work is motivated by some beautiful conjectures/open problems
in the field. It is well-known that any {\qfm} admits at least one closed, embedded, and incompressible {\ms}.
The following question, probably due to Hass-Thurston (see \cite{GW07}) and Uhlenbeck \cite{Uhl83}), addresses the multiplicity question:

\begin{conj}\label{HT}
For any integer $N> 0$, and any closed surface $S$ of genus $g\ge 2$,
does there exist a {\qfg} $G\cong\pi_{1}(S)$ such that the
resulting {\qfm} $\Mscr = \H^3/G$ contains at least $N$ distinct,
immersed, closed, incompressible {\ms}s, all diffeomorphic to $S$?
\end{conj}

Note that Anderson (\cite{And83}) constructed a {\qfm}
containing at least two incompressible {\ms}s, and we
(\cite{HW15b}) have constructed, given any prescribed
positive integer $N$, a {\qf} manifold (\emph{whose genus depends on $N$}) containing at least
$N$ distinct, embedded, closed, incompressible, (locally {\la}) {\ms}s.

Ever since Thurston's geometrization theorem for fibered three-manifolds (see for instance \cite{Thu82}), their
geometry is extremely important in {\htm} theory. We will also investigate closed {\ms}s in closed {\htm}s that fiber over the circle.

\begin{Def}\label{mf}
We call a $C^2$-fibration \textbf{minimal} or \textbf{geometrically taut} on an oriented closed {\htm} $\Mcal$
that fibers over the circle with fiber $S$ if each leaf is a closed incompressible {\ms}, which is
homeomorphic to the fiber $S$.
\end{Def}

By a celebrated theorem of Sullivan (\cite{Sul79}), any closed
Riemannian manifold with \emph{taut foliation} (a codimensional
one $C^2$-foliation such that there is a closed loop transversal
to each leaf) admits a minimal foliation {\wrt} some Riemannian
metric. The existence of a {\mf} structure has tremendous applications
in Riemannian geometry. A famous question in this direction is the following (see for instance \cite{Has05, Rub07}):

\begin{conj}\label{Th}
Does there exist a {\fhtm} which admits a minimal foliation?
\end{conj}

These questions are intricately related, see for instance \cite{And83}, where Anderson further conjectured that any closed
{\htm} does not admit a local parameter family of closed {\ms}s, in particular, does not admit a foliation of closed {\ms}s. These
questions have had profound impact in the theory of {\htm}s, as well as many other fields. In this work, we address problems
related to these questions.

\subsection{Main results}
In this paper, we analyze the relationship between the {\cl} of {\scg}s in a metrically complete {\htm} (essentially just inside solid
tubes) and closed {\la} {\ms}s in such {\htm}s. In one dimension lower, when a {\scg} $\gamma$ is short enough, any closed
geodesic disjoint from $\gamma$ can not go too deep inside the collar neighborhood of $\gamma$. Intuitively our argument is
similar in spirit, but we need to involve the {\cl} (real length and twist angle) to prevent a closed incompressible {\la} {\ms} going
too deep into a {\mst}. As a consequence of this relationship, we prove statements regarding the multiplicity of closed
incompressible {\ms}s in {\qfm}s, and the (non)existence of minimal fibrations on certain oriented closed {\htm}s that
fiber the circle.

Before we state our main results, we define some constants that will appear in the main statements which play an essential role
in our argument. These constants are unified through the following function:

\begin{Def}
We define the function $\Ot(x): [1,\infty) \to (0,1)$ as follows:
\be\label{W}
\Ot(x) =\frac{\sqrt{3}}{4\pi}
          \left[\cosh^{-1}\left(\frac{1}{1+
          \sqrt{1+(8x^2-8x+1)^2}}+1\right)\right]^2\ .
\ene
\end{Def}

It is elementary to verify that $\Ot(x)$ is a decreasing function of
$x \in [1,\infty)$, and $\lim\limits_{x\to \infty}\Ot(x) = 0$.
The maximum value is $\Ot(1) \approx 0.107071$, a fundamental constant
in {\htm} theory: Meyerhoff's constant.

Now we define the following ``Otal's constant", depending only
on the genus $g \ge 2$ of $S$:

\be\label{Otal}
\resizebox{0.875\hsize}{!}{$
   \varepsilon_{\Otal}(g)
      = \Ot(g)
      = \frac{\sqrt{3}}{4\pi}
          \left[\cosh^{-1}\left(\frac{1}{1+
          \sqrt{1+(8g^2-8g+1)^2}}+1\right)\right]^2\ .
   $}
\ene

Otal (\cite{Ota95, Ota03}) showed that when a curve (i.e., {\scg})
is sufficiently short, it is unknotted in a natural sense, and we always have
$0 < \varepsilon_{\Otal}(g) \leq \varepsilon_{\Otal}(2) =\Ot(2) \approx{}0.01515$.
We prove the following theorem on the multiplicity of closed
{\ms}s in {\qf} manifolds:

\bt\label{main1}
If an oriented closed {\htm} $\Mcal$ that fibers over the circle
with fiber $S$ contains a {\scg} whose {\cl}
$\Lcal = \ell + \sqrt{-1}\,\theta$ satisfies:
\ben
\item $\ell < \varepsilon_{\Otal}(g)$;
\item \be\label{ratio}
         \frac{|\theta|}{\sqrt{\ell}} >
         \sqrt[4]{3\pi^2} \approx 2.33268\ ,
      \ene
\een
then for any positive integer $N$, there exists a {\qfm}
$\Mscr \cong S\times \R$ which contains
at least $N$ embedded closed incompressible {\la} {\ms}.
\et

The techniques developed in \cite{HW15b} do not extend to the case of
arbitrary genus. Theorem~\ref{main1} states that for \emph{ANY} genus
$g \geq 2$, assuming the above two conditions on the {\cl} of some short
curve on an oriented closed {\htm} $\Mcal$ that fibers over the circle
with fiber $S$, then one can find a {\qfm} $\Mscr$ which contains arbitrarily
many  embedded closed incompressible {\ms}s. For different integer $N$'s,
the {\qfm}s obtained from this scheme are possibly different. This result
is also an improvement from \cite{Wan12}. Via computer programs, in
Appendix \S5.2, explicit examples of {\fhtm}s which satisfy the conditions in
Theorem ~\ref{main1} are produced. Higher genera cases follow from standard
finite coverings of lower genera examples. Therefore, we obtain an affirmative
answer to Question ~\ref{HT}.

It is well-known (\cite{Thu80} or \cite[Corollary 4.3]{Thu98}) that
$\lim\limits_{\ell\to 0}\theta = 0$ (we provide a proof in the
Appendix of this paper), but their quantitative nature for short
curves is notoriously difficult to control.
Minsky (\cite[Lemma 6.4]{Min99}) obtained a uniform upper bound
for any {\scg} in a Kleinian surface group with {\cl}
$\Lcal = \ell + \sqrt{-1}\,\theta$:

\be\label{minsky}
\frac{|\theta|}{\sqrt{\ell}} < \sqrt{\frac{2\pi}{C_1}}\ ,
\ene
where $C_1$ is a positive constant depending only on $g$. See
Corollary ~\ref{rot} in the appendix for an explicit bound.

Next we define a universal constant
\be\label{qf}
   \varepsilon_{0}=\frac{\sqrt{3}}{4\pi}
          \left[\cosh^{-1}\left(\frac{1}{1+
          \sqrt{1+(7+4\sqrt{3})^2}}+1\right)\right]^2
          \approx{}0.01822\ .
\ene
Using above $\Ot$-function notation \eqref{W}, one verifies that
$\varepsilon_0 = \Ot\left(\frac{2+\sqrt{3}}{2}\right)$ and
$\Ot(g) < \varepsilon_0$ for all $g \ge 2$.

In this paper, whenever we mention a foliation or fibration on a
mapping torus, we always assume it is $C^2$, and each leaf is a
closed surface diffeomorphic to the surface $S$ which is used to defined
the mapping torus, as we apply results of Sullivan (\cite{Sul79}),
Harvey-Lawson (\cite{HL82}) and Hass (\cite{Has86}) in an essential way
for our next result. We prove the following result related to
Question ~\ref{Th}:

\bt\label{main2}
If an oriented closed {\htm} $\Mcal$ that fibers over the circle
with fiber $S$ contains a {\scg} whose {\cl}
$\Lcal = \ell + \sqrt{-1}\,\theta$ satisfies:
\ben
\item $\ell < \varepsilon_{\Otal}(g)$;
\item $|\theta|/\sqrt{\ell} > \sqrt[4]{3\pi^2} \approx 2.33268$,
\een
then $\Mcal$ does not admit a {\mf}.
\et

Explicit examples of closed {\fhtm}s fibering over the circle which satisfy our conditions in Theorem ~\ref{main2} are produced also
in \S5.2. As an immediate corollary, we have:
\begin{cor}\label{corThu}
There exists some {\fhtm}s which do not admit any minimal foliation (in the sense of Definition ~\ref{mf}).
\end{cor}

\br
Recently Hass (\cite{Has15}) also obtained results on the question ~\ref{Th} using arbitrarily short geodesics on {\fhtm}s. We are
thankful for the correspondence.
\er
\subsection{Comments on the techniques and constants}
Margulis tubes are fundamental tools in {\tm} theory,
but it is usually very difficult to carry out explicit
calculations using {\mt}s of short curves in the study of {\htm}s.
We work with {\emph {\mst}s} (see \cite{Mey87}) instead in this
paper since we seek more computable conditions.

Otal's constant $\varepsilon_{\Otal}(g) = \Ot(g)$ did not
directly appear in his work \cite{Ota95}. In order to show
a sufficiently short geodesic $\gamma$ is unknotted, he
requires that the area of the meridian disk of the {\mt} of
$\gamma$ is greater than $4\pi(g-1)$. In our argument,
we replace the role of {\mt} by the \emph{\mst} of $\gamma$,
and we require, if $\ell$ (the real length of $\gamma$)
is less than this ``Otal's constant", then the area of the
meridian disk of the {\mst} of
$\gamma$ is greater than $4\pi(g-1)$ (See
Proposition ~\ref{meridian}). The numerical number is
calculated following this idea and using Meyerhoff's constant.

The other constant
$\varepsilon_0 = \Ot\left(\frac{2+\sqrt{3}}{2}\right)$ is designed so that a
{\la} {\ms} constructed similar to Calegari-Gabai (\cite{CG06})
by the means of shrink-wrapping will be
separated from the core curve in the {\mst} $\T$ whose {\cl}
satisfies the conditions in Theorem ~\ref{main3} (see
Lemma ~\ref{claim}), a main ingredient in the proof of both main Theorems ~\ref{main1} and ~\ref{main2}.

There are two other constants that will appear later. One is
$\varepsilon_1$ in the statement of Theorem
~\ref{area-comp}. Using our $\Ot$-function in \eqref{W},
we note here $\varepsilon_1 = \Ot\left(\frac32\right)\approx 0.03347$.
The other is Meyerhoff's constant $\varepsilon_2 = \Ot(1)$
which appears in Theorem ~\ref{Mey}. This is to
guarantee the existence of the {\mst}s around short curves.

In terms of our technical needs, we need $\ell < \varepsilon_2$ to define {\mst}s for short curves, and we
need a stricter $\ell < \varepsilon_1$ for a technical reason in a key inequality \eqref{area-comp2} in Theorem
~\ref{area-comp}. We need the above mentioned separation between a closed {\ms} and a short curve,
established using an even stricter condition $\ell < \varepsilon_0$, to prove Theorem ~\ref{main3}, and finally
we require further $\ell < \varepsilon_{\Otal}(g)$ in the proof of Theorems ~\ref{main1} and ~\ref{main2} to prevent curves from being
knotted or linked. In short, we have the following ordered constants which control the real length of a short geodesic:
\be\label{order}
\resizebox{0.875\hsize}{!}{$
   \varepsilon_{\Otal}(g)
   \leq \varepsilon_{\Otal}(2)
   < \varepsilon_0 = \Ot\left(\frac{2+\sqrt{3}}{2}\right)
   < \varepsilon_1 = \Ot\left(\frac32\right)
   < \varepsilon_2 = \Ot(1)\ .
$}
\ene
\subsection{Outline of the paper}
The organization of the paper is as follows: we summarize necessary background on {\ksg}s, mapping tori,
minimal helicoids in $\H^3$ and the {\mst} around short curves in {\htm}s in \S2. We develop our methods in
\S3, and use these techniques to prove our main theorems in \S4. We include an appendix to include a proof of a
Proposition by Thurston on {\cl} of short curves which provides an upper bound for the ratio $|\theta|/\sqrt{\ell}$.
Also in the appendix, we use Twister (\cite{Twister}) and SnapPy (\cite{SnapPy}) programs to produce some explicit
examples of fibered {\htm}s with our conditions satisfied. Taking finite covers of these lower genera examples one finds
more examples for higher genera surfaces of which mapping tori are made.

\subsection{Acknowledgement}
We are grateful to Saul Schleimer for his generous help on producing some explicit examples using Twister program. We also wish to thank Mark Bell for
his help with Twister. B. W. thanks the support from PSC-CUNY Research Award \#{}68119-0046. Z. H. acknowledges supports from
U.S. NSF grants DMS 1107452, 1107263, 1107367 ``RNMS: Geometric Structures and Representation varieties" (the GEAR Network).
This work was supported by a grant from the Simons Foundation (\#359635, Zheng Huang). We thank ICERM at Brown University
where part of this work was conducted while participating the special program on ``Low dimensional geometry, topology and dynamics". Finally
we are very grateful to the referee for many detailed suggestions and corrections.

\section{Tool Box}
\subsection{{\ksg}s and hyperbolic mapping tori}
We will mostly work with the {\uhs} model of hyperbolic three-space:
$\H^3 = \{z + tj: z \in \C,\, t>0\}$,
equipped with the standard {\hym}:
$ds^2 = \frac{|dz|^2+dt^2}{t^2}$.
The orientation preserving isometry group of $\H^3$, denoted by
$\PSL$, is the set of M\"obius transformations on $\H^3$,
namely, for each element $\tau\in\PSL$, we have
\beq
\tau(z) = \frac{az+b}{cz+d}\ ,\ \ \  \forall z \in \C\ ,
\eeq
with $ad-bc=1$. Its Poincar\'e extension is given by:
\beq
\tau(z+tj) = \frac{(az+b)(\overline{cz+d})+a\bar{c}t^2+tj}
{|cz+d|^2+|c|^2t^2}\ ,\ \ \  \forall (z,t) \in \H^3\ .
\eeq

Suppose that $S$ is an oriented closed surface of genus $\geq{}2$.
Let $\rho:\pi_1(S)\to\PSL$ be a discrete and faithful representation,
then the image $G=\rho(\pi_1(S))$, a discrete subgroup of $\PSL$, is
called a \textbf{{\ksg}}. The quotient manifold
$M_\rho=\H^3/\rho(\pi_1(S))$ is a complete {\htm}.
By the work of Thurston and Bonahon (\cite{Bon86}), we know that
$M_\rho$ is diffeomorphic to $S\times\R$.

Two {\ksg}s are \textbf{equivalent} if the corresponding representations are conjugate in $\PSL$. The \textbf{\ads} of $S$,
denoted by $\AH(S)$, is the space of equivalence classes. A {\ksg} is called \textbf{\qf} if its limit set is a topological
circle. The resulting quotient of $\H^3$ by a {\qf} group is called a \textbf{\qfm}. We
abuse our notation to denote both the space of {\qfm}s and
the space of {\qf} groups by ${\QF}(S)$. This space
plays a fundamental role in {\htm} theory.


Let $\{\rho_n:\pi_1(S)\to\PSL\}$ be a sequence of representations, then the sequence of {\ksg}s
$\{G_n=\rho_n(\pi_1(S))\}$ \textbf{converges algebraically} if $\lim\limits_{n\to\infty}\rho_n(\gamma)$ exists as a
M\"obius transformation for all $\gamma \in \pi_1(S)$.
Since the Kleinian surface group is finitely generated, the algebraic
limit of {\ksg}s is also Kleinian (see \cite{JK82}). Equipping the
deformation space $\AH(S)$ with the topology of algebraic convergence, the
space $\AH(S)$ is closed as a subspace of the space of equivalence classes of all homomorphisms into the isometry
group of $\H^3$ (see \cite{Chu68,Wie77} or \cite[Proposition 1.1.3]{Ota01}). One of the
fundamental theorems in {\ksg} theory is that ${\QF}(S)$ is in fact the interior of $\AH(S)$ (see \cite{Mar74, Sul85, Min03}).
Moreover, if we denote $\overline{\QF(S)}$ the closure of
${\QF}(S)$ in $\AH(S)$ {\wrt} the algebraic topology,
then we have (see \cite{BB04, Bro07}):
\be\label{closure}
\overline{\QF(S)} = \AH(S)\ .
\ene

A mapping torus with monodromy $\psi:S\to{}S$, denoted by $\Mcal_\psi$, can be constructed
by taking the quotient $S\times[0,1]/\sim$, where we identify $(x,0)$ and $(\psi(x),1)$. The automorphism $\psi$ of
$S$ defines an element of the mapping class group $\Mod$, it is \textbf{pseudo-Anosov} if no power of $\psi$ preserves the
isotopy class of any essential {\scg} on $S$. Thurston's hyperbolization theorem
(see \cite{Thu98,Mcm96,Ota96,Ota01})
shows that the mapping torus $\Mcal_\psi$ carries a hyperbolic structure if and only if $\psi$ is pseudo-Anosov, in this case
$\Mcal_{\psi}$ or simply $\Mcal$, is an oriented closed {\htm} that fibers over the circle with fiber $S$.
Though the hyperbolic mapping tori and {\qfm}s are very different geometrically, Thurston has shown that a certain covering of
the hyperbolic mapping tori arises as the limit of {\qfm}s:
Let $\Mcal_\infty$ be the infinite
cyclic cover of $\Mcal$ corresponding to the subgroup $\pi_1(S) \subset \pi_1(\Mcal)$, then $\Mcal_\infty$ is a doubly degenerated {\htm}
diffeomorphic to $S \times \R$ which arises as an algebraic limit of manifolds in  ${\QF}(S)$, hence it lies on the boundary of {\qf} space within $\AH(S)$.
\subsection{Family of Helicoids in hyperbolic three-space}

First let us describe a construction of a helicoid in $\H^3$,
which will descend to a {\ma} in a {\mst} in \S{}3.2.

\begin{Def}
The \textbf{helicoid} $\Hscr_a$ in $\H^3$, the {\uhs} model of the hyperbolic $3$-space, is the surface parametrized
by the $(u,v)$-plane as follows:
\be\label{heli}
\Hscr_a = \left\{z+tj \in \H^3: z = e^{v+\sqrt{-1}\,av}\tanh(u)\ ,\
t=\frac{e^v}{\cosh(u)} \right\}\ ,
\ene
where $-\infty<u,v<\infty$. In this model, the axis of $\Hscr_a$ is the $t$-axis.
\end{Def}

The first fundamental form can be written as
\be\label{heli-I}
I = du^2 + (\cosh^2(u)+a^2\sinh^2(u))dv^2\ .
\ene
Mori proved $\Hscr_a$ is indeed a {\ms} in $\H^3$ (see \cite{Mor82}). We will soon choose
$a = |\theta|/\ell$. We want to remark that this {\ms} is a beautiful analog of the helicoid
in {\eus}, namely, it is a ruled surface (see \cite{Tuz93}) that is stratified into straight lines {\wrt}
the {\hym}. This property is used in the proof of our key Lemma ~\ref{claim}.

\begin{figure}[htbp]
\centering
\begin{minipage}[t]{1\linewidth}
\centering
\includegraphics[scale=0.2]{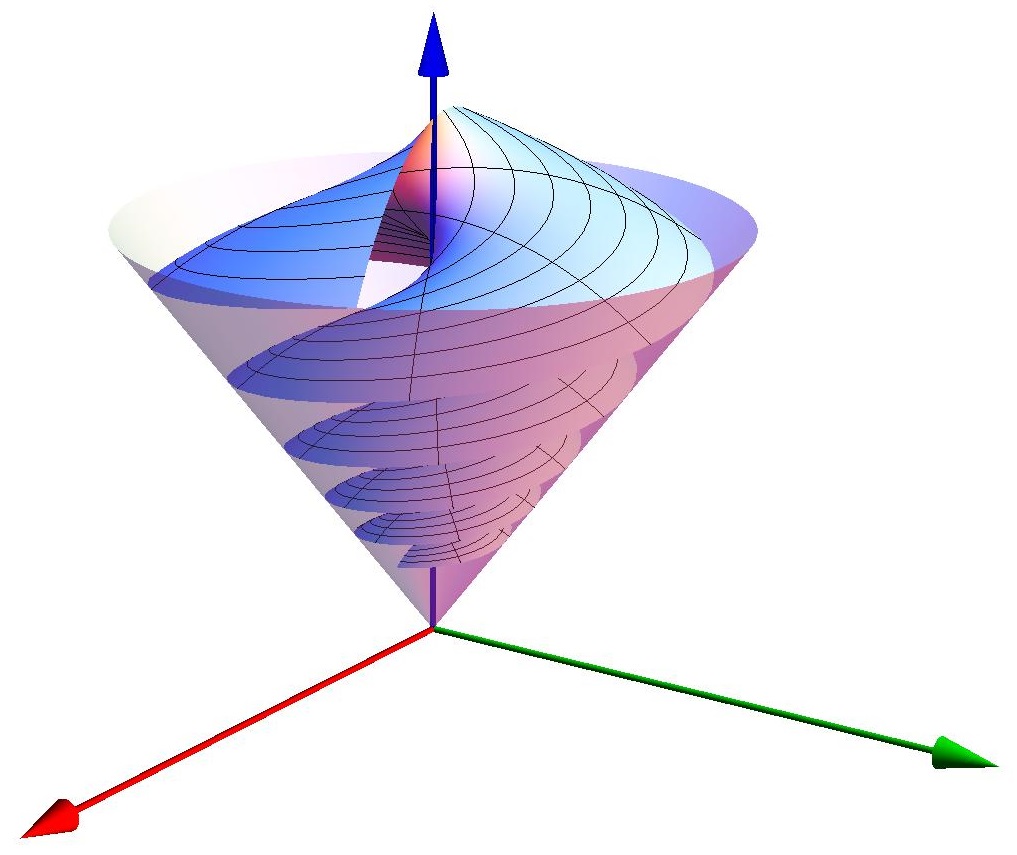}
\end{minipage}%
\caption{The helicoid $\Hscr_{10}$ defined by \eqref{heli}
for $-\log(2)\leq{}u\leq\log(2)$ and $0\leq{}v\leq\log(5)$
in the {\uhs} model. The cone is the $\log(2)$-neighborhood of
the $t$-axis. The curves perpendicular to the spirals are
geodesics in $\H^3$.}
\end{figure}

\section{Tubes of short curves and minimal surfaces}
In this section, we start to develop techniques for proving our main theorems. We will work with {\mst}s
associated with short {\scg}s in complete {\htm}s.

\subsection{Short curves and deep tubes}

As mentioned in the introduction, we make use of the {\mst}s
around short curves, instead of the
{\mt}s. This approach makes our calculations more explicit.
In this subsection, we construct such
{\mst}s (following \cite{Mey87}).

We consider loxodromic elements in the Kleinian
surface group, namely, $\tau(z) = \alpha z$, up to conjugacy,
where $\alpha = \exp(\ell+\sqrt{-1}\,\theta)$ with $\ell > 0$
and $\theta \in [-\pi,\pi)$. Such a loxodromic
element translates points on the $t$-axis by the (hyperbolic)
distance $\ell$ and twists a normal plane
by the angle $\theta$. For a {\scg} $\gamma$ in any complete
{\htm}, a lift $\tilde{\gamma}$ of $\gamma$ in $\H^3$
is the axis of a loxodromic element $\tau \in \PSL$ representing
$\gamma$. We note that a different lift
just gives rise to another element that is conjugate to $\tau$
in $\PSL$.

\begin{Def}\label{mst}
We denote
\beq
  \Nscr_{r}(\tilde{\gamma}) = \{x \in \H^3:
  \dist(x,\tilde{\gamma}) < r\}\ ,
\eeq
as the $r$-neighborhood of the geodesic $\tilde{\gamma}$ in $\H^3$.
We call $r_0(\gamma)$ the \textbf{tube radius} of $\gamma$ if it is the supremum of numbers $r>0$ such that $\Nscr_r(\tilde{\gamma}) \cap \Nscr_r(\tilde{\gamma}') =\emptyset$,
for all lifts $\tilde{\gamma}'$ of $\gamma$ different from
$\tilde{\gamma}$. The \textbf{\mst} of $\gamma$
is then defined by, for $\tau$ loxodromic in $\PSL$ representing $\gamma$, and $\tilde{\gamma}$ the axis of the element $\tau$,
\be\label{mst2}
   \T= \Nscr_{r_0}(\tilde{\gamma})/\langle\tau\rangle\ .
\ene
\end{Def}

We have the following basic result of Meyerhoff:

\bt[{\cite{Mey87}}]\label{Mey}
If $\gamma$ is a {\scg} in a complete {\htm} with real length
$\ell$ less than the constant
\be\label{MeyCon}
   \varepsilon_2 = \Ot(1) =
   \frac{\sqrt{3}}{4\pi}\left(\log(\sqrt{2}+1)\right)^2 \approx 0.107071\ ,
\ene
then there exists an embedded {\mst} around $\gamma$ whose
tube radius is given by
\be\label{tuberad}
   \cosh^2(r_0(\gamma)) = \frac12
   \left(\frac{\sqrt{1-2\kappa(\ell)}}{\kappa(\ell)}+1\right)\ ,
\ene
where the function
\be\label{kappa}
   \kappa(\ell) =\cosh\left(\sqrt{\frac{4\pi\ell}
   {\sqrt{3}}}\right) - 1\ .
\ene
Moreover, {\mst}s around different {\scg}s do not intersect if
their real lengths are both less than $\varepsilon_2$.
\et
Note that Meyerhoff's constant \eqref{MeyCon} is the maximum of
the $\Ot$-function, therefore the real
length condition in Theorem ~\ref{Mey} is satisfied by short
curves in our main results.

We now justify the geometry behind the introduction of
``Otal's constant":

\bpro\label{meridian}
Let $\gamma$ be a {\scg} in a complete {\htm} which is
diffeomorphic to $S \times\R$, such that its real length
$\ell$ is less than ``Otal's constant" \eqref{Otal}, namely,
$\ell < \varepsilon_{\Otal}(g) =\Ot(g)$, where $g \ge 2$
is the genus of $S$, then the area of the meridian disk of the
{\mst} $\T$, defined in \eqref{mst2}, is greater than
the hyperbolic area of $S$.
\epro

\bp
Recall that
\beq
   \varepsilon_{\Otal}(g)= \Ot(g) = \frac{\sqrt{3}}{4\pi}
          \left[\cosh^{-1}\left(\frac{1}{1+
          \sqrt{1+(8g^2-8g+1)^2}}+1\right)\right]^2\ .
\eeq
From \eqref{kappa}, we have
\beq
   \kappa(\ell) < \frac{1}{1+ \sqrt{1+(8g^2-8g+1)^2}}\ .
\eeq
Then the tube radius $r_0(\gamma)$ satisfies
\beq
   \cosh^2(r_0(\gamma))=\frac12
   \left(\frac{\sqrt{1-2\kappa(\ell)}}{\kappa(\ell)}+1\right)
   >\frac{8g^2-8g + 2}{2}\ .
\eeq
The area of the meridian disk of the tube $\T$ is
$2\pi(\cosh(r_0)-1)$. Therefore we have
\be\label{eq:Otal}
   2\pi(\cosh(r_0)-1) > 2\pi(2g-2)\ ,
\ene
which is the hyperbolic area of the surface $S$.
\ep

With this estimate, if $\ell < \varepsilon_{\Otal}(g) =\Ot(g) < \Ot(1)$,
the arguments in \cite{Ota95, Ota03} imply that $\gamma$ is unknotted
(i.e., $\gamma$ lies on an embedded surface isotopic to $S$), and any
collection of {\scg}s with real lengths shorter than $\Ot(g)$ are not
linked (i.e., they can be isotoped such that each lies on a surface
$S \times n$, where $n$ is an integer).

Using the computation in
Proposition \ref{meridian}, we also have the following proposition, which
will be used for proving Theorem \ref{main1} and Theorem \ref{main2}.

\bpro\label{pro:unknotted}
Let $\psi:S\to{}S$  be a pseudo-Anosov map of a closed surface $S$
with genus $g\geq{}2$, and $\Mcal_{\psi}$ be the hyperbolic mapping
torus with respect to $\psi$.
If $\Mcal_{\psi}$ contains a {\scg} $\gamma$ whose real length satisfies
$\ell < \varepsilon_{\Otal}(g)$, then $\gamma$ lies on an embedded
surface isotopic to $S$.
\epro

\bp
Recall that the fundamental group of $\Mcal_{\psi}$ is given by
\begin{equation}\label{eq:fundamental_group_mapping_torus}
   \pi_{1}(\Mcal_{\psi})=\langle{}\pi_{1}(S),t\ :\
   \forall\alpha\in\pi_{1}(S),\ t\alpha{}t^{-1}=\psi_{*}(\alpha)\rangle\ ,
\end{equation}
where $\psi_{*}:\pi_{1}(S)\to\pi_{1}(S)$ is the induced isomorphism.

According to \cite{SY79, SU82}, $\Mcal_{\psi}$ contains an embedded least area
incompressible minimal surface $\Sigma$, which is isotopic to some fiber
of $\Mcal_{\psi}$. By the Gauss-Codazzi equation, we have
$K_{\Sigma}=-1-\lambda^{2}$, where $K_{\Sigma}$ and $\pm\lambda$ are the
Gauss curvature and principal curvatures of $\Sigma$ with respect to the induced
metric respectively.
By Gauss-Bonnet Theorem, the area of $\Sigma$
is bounded from above by the area of a hyperbolic closed surface of genus $g$, i.e.,
$\Area(\Sigma)\leq{}2\pi(2g-2)=4\pi(g-1)$.

Obviously the minimal surface $\Sigma$ is non-separating in $\Mcal_{\psi}$,
so there exists an essential loop $\gamma'\subset\Mcal_{\psi}$ intersects
$\Sigma$ transversely at a single point. Let $\gamma^{*}\subset\Mcal_{\psi}$
be the simple closed geodesic homotopic to $\gamma'$, then the mod $2$
intersection number of $\gamma^{*}$ and $\Sigma$, denoted by
$I_{2}(\gamma^{*},\Sigma)$, is equal to one.

We claim that the length of $\gamma^{*}$ is $\geq\varepsilon_{\Otal}(g)$.
Our argument below actually can show that any such kind of geodesic has
length greater than $\varepsilon_{\Otal}(g)$. In particular, the geodesic
representing the element $t$ in \eqref{eq:fundamental_group_mapping_torus}
has length greater than $\varepsilon_{\Otal}(g)$.

Let $\mathbb{T}(\gamma^{*})\subset\Mcal_{\psi}$ be the {\mst} of $\gamma^{*}$
with tube radius $r_{0}$. Similar to Lemma \ref{lem:Lemma_5_2_HW15a},
any component of $\Sigma\cap\Nscr_{s}(\gamma^{*})$ is either a disk or an
annulus for $0\leq{}s\leq{}r_{0}$. Recall $I_{2}(\gamma^{*},\Sigma)=1$, so
for any $0\leq{}s\leq{}r_{0}$, $\Sigma\cap\Nscr_{s}(\gamma^{*})$ must contain
at least one disk component, whose boundary is a simple closed curve
which is essential on $\partial\Nscr_{s}(\gamma^{*})$ and is
isotopic to the meridian of $\partial\Nscr_{s}(\gamma^{*})$ .
The meridian of $\partial\Nscr_{s}(\gamma^{*})$ is a geodesic with respect to
the induced metric on the torus $\partial\Nscr_{s}(\gamma^{*})$,
whose length equals $2\pi\sinh(s)$, so we have the following inequality
\begin{equation*}
   \Length(\Sigma\cap\partial\Nscr_{s}(\gamma^{*}))\geq{}2\pi\sinh(s)\ ,
   \ 0\leq{}s\leq{}r_{0}\ .
\end{equation*}
By the coarea formula \cite[p.399]{CG06}
or \cite[Lemma 5.4]{Wan16}, we have
\beq
\begin{aligned}
   \Area(\Sigma\cap\mathbb{T}(\gamma^{*}))
     & = \int_{0}^{r_{0}}\int_{\Sigma\cap\partial\Nscr_{s}(\gamma^{*})}
         \frac{1}{\cos\theta}dlds \\
     & \geq\int_{0}^{r_{0}}\Length(\Sigma\cap\partial\Nscr_{s}(\gamma^{*}))ds \\
     & \geq\int_{0}^{r_{0}}2\pi\sinh(s)ds
       = 2\pi(\cosh(r_{0})-1)\ ,
\end{aligned}
\eeq
where $\theta=\theta(q)$ is the angle between the tangent space to $\Sigma$ at $q$,
and the radial geodesic that is through $q$ (emanating from $\gamma^{*}$) and is
perpendicular to $\gamma^{*}$.
If the real length of $\gamma^{*}$ satisfies $\ell<\varepsilon_{\Otal}(g)$,
then by the computation in Proposition \ref{meridian}, we have the inequality
\eqref{eq:Otal}, which implies the inequality
$\Area(\Sigma\cap\mathbb{T}(\gamma^{*}))>4\pi(g-1)$.
But this is a contradiction, since we have $\Area(\Sigma\cap\mathbb{T}(\gamma^{*}))\leq\Area(\Sigma)\leq{}4\pi(g-1)$.

Therefore if $\gamma\subset\Mcal_{\psi}$ is a geodesic with
$\ell < \varepsilon_{\Otal}(g)$, then the mod $2$ intersection number of
$\gamma$ and $\Sigma$ (or any fiber of $\Mcal_{\psi}$) is zero, so we may
lift $\gamma$ to a closed geodesic, denoted by $\widetilde\gamma$, in the 
cyclic covering space of $\Mcal_{\psi}$ which is
diffeomorphic to $S\times\R$.
Now we can apply the results in \cite{Ota95,Ota03}
to show that $\widetilde\gamma$ lies on an embedded surface isotopic to $S\times\{0\}$.
In the quotient space $\Mcal_{\psi}$, therefore $\gamma$ is contained in some fiber
of $\Mcal_{\psi}$.
\ep

\subsection{Minimal annuli in {\mst}s}

We now construct a minimal annulus inside a {\mst} of a short {\scg} in a metrically complete {\htm}, this is done
by using the helicoid in $\H^3$ defined in \eqref{heli}.

Let $\gamma$ be a {\scg} in a complete hyperbolic $3$-manifold $M$, and $\T$ be its {\mst} with tube
radius $r_0$. Let $\tau\in\PSL$ be a loxodromic element of {\cl} $\Lcal = \ell + \sqrt{-1}\,\theta$
representing $\gamma$, with $\ell > 0$ and $\theta \in [-\pi,\pi)$. Suppose $\tilde{\gamma}$ is a lift of
$\gamma$ in $\H^3$ which is the axis of $\tau$, letting $a = |\theta|/\ell$, we define a surface in $\T$ as follows:
\be\label{annulus}
   \Ascr_a =\frac{\Hscr_a\cap\Nscr_{r_0}(\tilde{\gamma})}
   {\langle \tau \rangle}\ .
\ene
It is not hard to see that each component of
$\Ascr_a \cap\partial \Nscr_r(\gamma)$ is a closed geodesic
{\wrt} the induced metric on $\partial \Nscr_r(\gamma)$,
for each $r \in (0, r_0]$, with
$\Nscr_r(\gamma) = \Nscr_r(\tilde{\gamma})/\langle \tau \rangle$.
It is proven in \cite{Wan12} that
$\Ascr_a$ is indeed a minimal annulus in $\T$, moreover,
its area is explicitly computed as:
\be\label{area}
   \Area(\Ascr_a) = 2\int_0^{r_0}\sqrt{\ell^2\cosh^2(u) +
   \theta^2\sinh^2(u)}\,du\ .
\ene

\begin{figure}[htbp]
\centering
\begin{minipage}[t]{1\linewidth}
\centering
\includegraphics[scale=0.2]{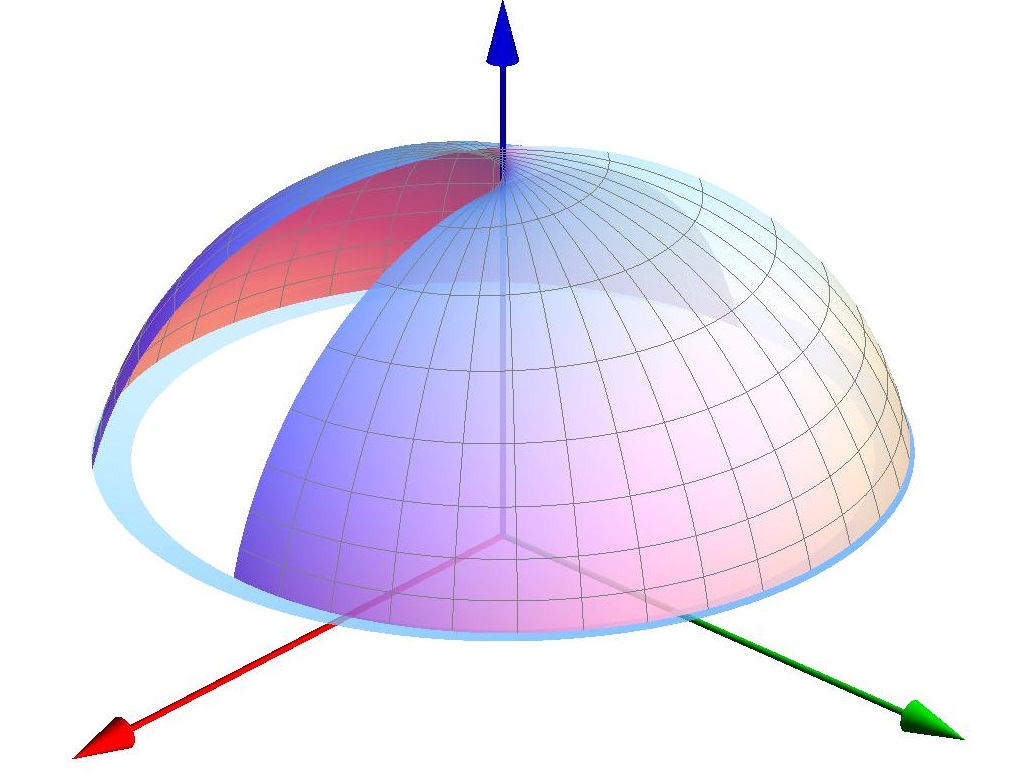}
\end{minipage}%
\caption{The fundamental domains of the maximal tube
of the closed geodesic $\gamma$ with complex length
$\Lcal=0.01+0.25i$ (the radius $r_0\approx{}1.98272$)
and the lifting of the minimal annulus $\Ascr_{25}$
contained in a piece of the helicoid $\Hscr_{25}$
which is given by \eqref{heli-I} for $a=25$,
$-r_0\leq{}u\leq{}r_0$ and $0\leq{}v\leq{}0.1$.
In this case, $0.25/\sqrt{0.01}=2.5$,
$\Area(\partial\T)\approx{}0.828202$ and
$\Area(\Ascr_{25})\approx{}1.35306$.}
\end{figure}
\vskip 0.1in
Now we prove the following technical estimate, where we
introduce a constant $\varepsilon_1 = \Ot\left(\frac32\right)$ to
guarantee a key inequality \eqref{area-comp2}, when the real
length of $\gamma$ is less than this constant
and the inequality \eqref{ratio} is satisfied.

\bt\label{area-comp}
If a complete {\htm} $M$ contains a {\scg} $\gamma$ whose {\cl}
$\Lcal = \ell + \sqrt{-1}\,\theta$ satisfies:
\ben
\item $\ell < \varepsilon_1 = \Ot\left(\frac32\right) =
       \frac{\sqrt{3}}{4\pi}\left[\cosh^{-1}
       \left(\frac{1}{1+5\sqrt{2}}+1\right)\right]^2
       \approx 0.03347\ $;
\item $|\theta|/\sqrt{\ell} >\sqrt[4]{3\pi^2} \approx 2.33268\ $,
\een
then we have
\be\label{area-comp2}
 \Area(\partial\T) = \pi\ell\sinh(2r_0) <
 |\theta|\cosh(r_0) < \Area(\Ascr_a)\ ,
\ene
where $a=|\theta|/\ell$.
\et

\bp
The area formula $\Area(\partial\T) = \pi\ell\sinh(2r_0)$
is well-known, see for instance \cite[Lemma 1.4]{GMM01}.

Recall from \eqref{kappa}, we have, once $\ell < \varepsilon_1$,
\be\label{kappa2}
   \kappa(\ell) =
   \cosh\left(\sqrt{\frac{4\pi\ell}{\sqrt{3}}}\right) - 1
   < \frac{1}{1+5\sqrt{2}}\ .
\ene
By the tube radius formula in \eqref{tuberad}, we then have:
\beq
\cosh(r_0) > \sqrt{\frac12\left(1+\left(1+5\sqrt{2}\,\right)
\sqrt{1-\frac{2}{1+5\sqrt{2}}}\,\right)} = 2\ .
\eeq
Applying the explicit area formula for the minimal annulus
\eqref{area}, we find:
\beqar
  \Area(\Ascr_a)
      &=& 2\int_0^{r_0}\sqrt{\ell^2\cosh^2(u) +
           \theta^2\sinh^2(u)}\,du \\
      &>& 2|\theta|\int_0^{r_0}\sinh(u)du
       = 2|\theta|(\cosh(r_0)-1) \\
      &\ge& |\theta|\cosh(r_0)\ .
\eeqar
Therefore to establish \eqref{area-comp2}, it suffices
to show
\be\label{lateruse}
   \ell^2\sinh^2(r_0) < \frac{\theta^2}{4\pi^2}\ .
\ene
First we note that $\kappa(\ell) \ge \frac{2\pi\ell}{\sqrt{3}}$.
Also $\ell < \varepsilon_1$ implies
\beq
   0<\sqrt{1-2\kappa(\ell)}-\kappa(\ell) < 1\ .
\eeq
Now we use \eqref{tuberad} again to find:
\beq
    \ell^2\sinh^2(r_0)
    =\frac{\ell^2}{2\kappa(\ell)}
     \left(\sqrt{1-2\kappa(\ell)}-\kappa(\ell)\right)
    <\frac{\ell^2}{2\kappa(\ell)}
    \le \frac{\sqrt{3}\ell}{4\pi}\ .
\eeq
We then complete the proof by using condition \eqref{ratio}.
\ep

\subsection{Separation} Let us recall some notations we will use in this subsection. Let $M=\H^3/G$ be a (metrically) complete
{\htm}, and let $\gamma$ be a {\scg} in $M$, whose real length
$\ell < \varepsilon_2$, the Meyerhoff constant. Also let $r_0$ be
the tube radius of $\gamma$, $\T$ be the {\mst} around $\gamma$, and
$\Nscr_r(\tilde{\gamma})$ (as in Definition ~\ref{mst}) be the
$r$-neighborhood of the lift $\tilde{\gamma}$ in $\H^3$, with
$r \in (0,r_0]$.

We will need the following result, proven in \cite{Wan12}, using
arguments similar to \cite{MY82, MY82b}, as well as \cite{FHS83}:

\bt[{\cite[Lemma 6]{Wan12}}]\label{disk}
Using the above notations, and let $\Nscr_r(\gamma) = \Nscr_r(\tilde{\gamma})/\langle\tau\rangle\subset{}M$,
where $\tau$ is the element in $G$ representing the geodesic $\gamma$.
If $C$ is a smooth simple curve which is null-homotopic on
$\partial \Nscr_r(\gamma)$ whose length is less than $2\pi\sinh(r)$,
with $0<r<r_0$, then $C$ bounds an embedded {\la} minimal disk
$\Delta \subset \Nscr_r(\gamma)\backslash\gamma$.
\et

For any $r \in (0,r_0]$, we let $D(r)$ be a disk on
$\partial \Nscr_r(\gamma)$ of radius $r$ in the induced metric on $\partial \Nscr_r(\gamma)$, and $\Bscr(r)$ be the {\la} minimal disk in $M$
bounding the closed curve $\partial D(r)$. We
define
\be\label{delta}
   \delta =\delta(\gamma,r_0) =\min\left\{\dist(\gamma, \Bscr(r)):
   \frac{r_0}{2} \le r \le r_0\right\},
\ene
where the distance is measured in {\hym}. We re-write
Theorem ~\ref{disk} into the following corollary to
quantify the separation of the minimal disk $\Delta$ and the
curve $\gamma$:

\bcor\label{coro}
Same notations as in above Theorem ~\ref{disk}. If $r \in \left[\frac{r_0}{2},r_0\right]$, then we have
$\delta \in \left(0,\frac{r_0}{2}\right)$, and the {\la}
disk $\Delta \subset \Nscr_r(\gamma)\backslash\gamma$, obtained
in Theorem ~\ref{disk}, satisfies $\dist(\gamma,\Delta) \ge \delta$.
\ecor


We prove the following existence result for a closed {\ms} with a specific property: it is separated from a {\scg}
if the curve satisfies our conditions in main theorems. More specifically,

\bt\label{main3}
Let $M$ be a closed or {\qf} {\htm}, and let $\gamma$ be a {\scg} contained in $M$ whose
{\cl} $\Lcal = \ell + \sqrt{-1}\,\theta$,
where $\ell>0$ and $\theta \in [-\pi,\pi)$, satisfies:
\ben
\item $\ell < \varepsilon_0 = \Ot\left(\frac{2+\sqrt{3}}{2}\right)\approx{}0.01822$;
\item $|\theta|/\sqrt{\ell}> \sqrt[4]{3\pi^2} \approx 2.33268$,
\een
if $S$ is an embedded closed incompressible surface of genus
$g\ge 2$ in $M\backslash\gamma$, then there exists an embedded {\la} {\ms}
$T\subset M\backslash\gamma$ isotopic to $S$. Here $T$ is of {\la}
means its area is the smallest among all surfaces in $M\backslash\gamma$ in the isotopy class of $S$.
\et

This can be easily applied to the case of {\qfm}s. Note that any {\qfm} always contains embedded, closed,
incompressible surfaces, a fact not always shared by some other classes of {\htm}s. The proof of this
theorem is along the lines of the arguments in \cite{Wan12}, but we need to take special care at places with
improved estimates. We also need the following lemma, whose proof can be found in
\cite{HW17}. We note that this is a purely topological lemma.

\bl[{\cite[Lemma 5.2]{HW17}}]\label{lem:Lemma_5_2_HW15a}
Let $\mathbf{T}$ be a solid torus with the core curve removed, and
let $S_{g,n}$ be a surface of type $(g,n)$ such that the Euler characteristic
of $S_{g,n}$ is negative
\footnote{The Euler characteristic of $S_{g,n}$ is $\chi(S_{g,n})=2-2g-n$.
Thus $\chi(S_{g,n})<0$ if one of the conditions is satisfied:
$g\geq{}2$ and $n\geq{}0$, $g=1$ and $n\geq{}1$ or $g=0$ and $n\geq{}3$.}.
If $S_{g,n}$ is embedded in
$\mathbf{T}$ such that $\partial{}S_{g,n}\subset\partial\mathbf{T}$, then there exists at least
one essential non-peripheral simple closed curve $\alpha\subset{}S_{g,n}$ such that
$\alpha$ bounds a disk $D\subset\mathbf{T}$.
\el

\bp[Proof of Theorem~\ref{main3}]
Our strategy will be first to invoke a technique modifying the {\hym} called ``shrink-wrapping", developed by
Calegari-Gabai (\cite{CG06}) in their work on the tameness conjecture. We use this to conformally modify the {\hym} of
$M$ inside a solid tube so that we can use the resulting {\tg} boundary tori as barriers. We then construct a {\ms} and
prove it is minimal {\wrt} the {\hym}.

Consider the solid torus $\Nscr_\sigma(\gamma) \subset M$, as
before, where $\sigma < r_0$ is a positive constant, and $r_0$
is the tube radius of $\gamma$. For each $t \in [0,1)$, one can
define a family of Riemannian metrics $g_t$ on $M$ such that it
coincides with the {\hym} on $M\backslash \Nscr_{(1-t)\sigma}(\gamma)$, while conformally
equivalent to the {\hym} on $\Nscr_{(1-t)\sigma}(\gamma)$.
Then by \cite[Lemma 1.18]{CG06}, for each $t \in [0,1)$, there is a function $f(t)$ satisfying
$\frac23 (1-t)\sigma< f(t) < \frac34 (1-t)\sigma$, such that the torus $\partial \Nscr_{f(t)}(\gamma)$ is {\tg} {\wrt} the metric
$g_t$, and the metric $g_t$ dominates the {\hym} on $2$-planes. Moreover, by the standard result in
\cite{SY79,FHS83,HS88}, for each $t \in [0,1)$, there exists an embedded surface $S_t$ in
$M\backslash \Nscr_{f(t)}(\gamma)$, isotopic to $S$, which is globally $g_t$-{\la} among all such surfaces. Our
first goal is to show, for $t$ sufficiently close to $1$, and $\gamma$ satisfying our conditions on {\cl}, the $S_t$
produced as a globally {\la} surface {\wrt} the metric $g_t$ is also of {\la} {\wrt} the {\hym}.

Letting $r \in [\frac{r_0}{2},r_0]$, by Corollary ~\ref{coro}, there is $\delta>0$ defined in \eqref{delta}, only
depending on $r_0$ and $\gamma$, such that the {\la} disk $\Delta \subset \Nscr_r(\gamma)\backslash\gamma$,
obtained in Theorem ~\ref{disk}, satisfies $\dist(\gamma,\Delta) \ge \delta$. Here $\delta \in (0,\frac{r_0}{2})$
and from above we have a constant $\sigma$ with $\sigma < r_0$. Now we choose $t$ sufficiently close to $1$
such that $(1-t)\sigma < \delta$. For instance we can define $t_0 \in (0,1)$, such that we have $(1-t)\sigma < \delta$ for all $t > t_0$.


Now we pause to prove a technical lemma which will be used later for our applications:

\bl\label{claim}
Under the conditions of Theorem ~\ref{main3}, for all $t > t_0$, where $t_0$ is the threshold defined above, the $g_t$-{\la} surface $S_t$
is disjoint from $\Nscr_\delta(\gamma)$, i.e.,
$S_t\cap \Nscr_\delta(\gamma) = \emptyset$, where $\delta$ is defined
in \eqref{delta}.
\el
\bp[Proof of Lemma ~\ref{claim}]
We may assume $S_t\cap \Nscr_{\delta}(\gamma) \not= \emptyset$,
for all $t > t_0$, and we will get a contradiction.

By the shrink-wrapping argument, it's easy to see that $S_{t}$ is disjoint from $\gamma$ for $t\in[0,1)$, hence
\begin{equation*}
   S_{t}\cap\T=S_{t}\cap\mathbf{T}\ ,
   \quad\text{for}\ t\in[0,1)\ ,
\end{equation*}
where $\T$ is the maximal solid tube of $\gamma$ with radius $r_0$
and $\mathbf{T}=\T\setminus\gamma$. Let $\Sigma$ be a component of $S_t\cap \T$ which intersects
$\Nscr_{\delta}(\gamma)$. We claim that we always have
\be\label{eq:large length}
   \Length(\Sigma\cap\partial \Nscr_s(\gamma))
   \ge 2|\theta|\sinh(s)\ ,
   \quad\text{where}\ \frac{r_0}{2}\le{}s\le{}r_0\ .
\ene
In fact, since $S_t$ is incompressible in $M$, and as we are dealing with a solid torus with core curve removed, we apply
Lemma \ref{lem:Lemma_5_2_HW15a}, then $\Sigma$ is either a disk whose boundary is null-homotopic on $\partial\T$ or an
annulus whose boundary is essential on $\partial\T$. There are two cases we need to consider:
\ben
\item \underline{Case One}: $\Sigma$ is a disk. Then by Corollary ~\ref{coro},
      $\Sigma\cap\Nscr_{s}(\gamma)$ consists of disjoint disks
      for all $s\in[\frac{r_0}{2},r_0]$. If there exists
      some $s'\in[\frac{r_0}{2},r_0]$ such that
      $\Length(\Sigma\cap\partial \Nscr_{s'}(\gamma))<2|\theta|\sinh(s')$,
      then $\Length(\Sigma\cap\partial \Nscr_{s'}(\gamma))<2\pi\sinh(s')$,
      which implies that $\Sigma\cap\Nscr_{s'}(\gamma)$ is disjoint from
      $\Nscr_{\delta}(\gamma)$, therefore so is $\Sigma$. A contradiction.
      Therefore \eqref{eq:large length} is true when $\Sigma$ is a disk.
\item \underline{Case Two}: $\Sigma$ is an annulus. For any
      $s\in[\frac{r_0}{2},r_0]$, $\Sigma\cap\Nscr_{s}(\gamma)$ either
      only consists of disjoint disks or contains an annulus, say $\Sigma'$.
      In the former subcase, \eqref{eq:large length} is true according to
      the argument in Case One. In the latter subcase, $\partial\Sigma'$
      consists of two isotopic slopes on $\partial\Nscr_{s}(\gamma)$, so
      $\Length(\partial\Sigma')$ is greater than either
      $2(2\pi\sinh(s))=4\pi\sinh(s)$ or
      $\Length(\Ascr_a\cap\partial \Nscr_s(\gamma))$,
      where $2\pi\sinh(s)$ is the length of the meridian of the torus
      $\partial\Nscr_{s}(\gamma)$ and each component of
      $\Ascr_a\cap\partial \Nscr_s(\gamma)$ is a geodesic isotopic to
      the longitude of the torus $\partial\Nscr_{s}(\gamma)$.
      It's easy to see the inequalities
      $4\pi\sinh(s)>2|\theta|\sinh(s)$ and
      \beqar
          \Length(\Ascr_a\cap\partial \Nscr_s(\gamma))
          &=& 2\sqrt{\ell^2\cosh^2(s)+\theta^2\sinh^2(s)}\\
          &>& 2|\theta|\sinh(s)\ .
      \eeqar
      So \eqref{eq:large length} is also true in this case.
\een

Recalling that the new metric $g_t$ dominates the {\hym} on
$2$-planes, we apply the co-area formula
(see \cite[Lemma 3]{Wan12}) to obtain the following estimate:
\bear\label{lowbd}
   \Area(\Sigma,g_t) &\ge& \Area(\Sigma)\notag\\
   &\ge&\Area\left(\Sigma\cap\left(\overline{\T\backslash
        \Nscr_{\frac{r_0}{2}}(\gamma)}\,\right)\right)\notag\\
   &\ge&\int_{\frac{r_0}{2}}^{r_0}
        \Length(\Sigma\cap\partial \Nscr_s(\gamma))ds\\
   &\ge&2|\theta|\int_{\frac{r_0}{2}}^{r_0}\sinh(s)ds \notag\\
     \label{lowbd2}
   &=& 2|\theta|\left(\cosh(r_0)-\cosh\left(\frac{r_0}{2}\right)\right)\ .
\eear
where we denote $\Area(\cdot,g_t)$ the $g_t$-area, and
$\Area(\cdot)$ the hyperbolic area.

We now interpret constant $\varepsilon_0$ in \eqref{qf}.
When $\ell < \varepsilon_0$, we have from \eqref{kappa} that:
\beq
   \kappa(\ell) =
   \cosh\left(\sqrt{\frac{4\pi\ell}{\sqrt{3}}}\right)-1< \frac{1}{1+\sqrt{1+(7+4\sqrt{3})^2}}\ ,
\eeq
therefore we have from Meyerhoff's formula \eqref{tuberad}
for the tube radius:
\be\label{eq:r0}
    \cosh(r_0) = \sqrt{\frac12 \left(\frac{\sqrt{1-2\kappa(\ell)}}
    {\kappa(\ell)}+1\right)} > \sqrt{3}+1\ .
\ene
As a result, we find:
\be\label{eq:cosh(r0)}
   \cosh(r_0) > 2\cosh\left(\frac{r_0}{2}\right)\ .
\ene
Putting this into the inequality \eqref{lowbd2}, we have:
\be\label{lowbd3}
   \Area(\Sigma)  > |\theta|\cosh(r_0)\ .
\ene
Since $\varepsilon_0 = \Ot\left(\frac{2+\sqrt{3}}{2}\right) <
\varepsilon_1 = \Ot\left(\frac32\right)$, conditions (i) and (ii) in the
statement allow us to apply the inequality \eqref{lateruse}
in the proof of Theorem ~\ref{area-comp}, namely,
we have
\beq
   \Area(\partial\T) = \pi\ell \sinh(2r_0) < |\theta|\cosh(r_0)\ .
\eeq
By our choice of $t> t_0$, we have $(1-t)\sigma < \delta$, the metric $g_t$ coincides with the {\hym} outside $\Nscr_{(1-t)\sigma}(\gamma)$,
and combining these inequalities, we have established:
\be\label{area-estimate}
   \Area(\Sigma,g_t)  > \Area(\partial\T)
   = \Area(\partial\T, g_t)\ .
\ene
This estimate then allows us to proceed with cut-and-paste, namely, we can replace each component of
$S_t\cap\T$ which intersects $\Nscr_{\delta}(\gamma)$ by either an annulus or a disk on $\partial\T$,
to obtain a new surface $S_t'\subset M\backslash \Nscr_\sigma(\gamma)$ such that it
has less $g_t$-area than $S_t$, disjoint from $\Nscr_\sigma(\gamma)$ and isotopic to $S$ in $M\backslash\gamma$.
This is impossible since $S_t$ is the {\la} surface with these properties. This completes the proof for the lemma.
\ep
Now we continue our proof for Theorem ~\ref{main3}. By above lemma,
we have the $g_t$-{\la} surface $S_t$ is separated from
$\Nscr_\delta(\gamma)$. But by shrinkwrapping, the metric $g_t$
coincides with the {\hym} outside of $\Nscr_\delta(\gamma)$, therefore
$S_t$ is of {\la} {\wrt} the {\hym} for $t$ sufficiently close
to $1$. Since $\delta$ is independent of $t$, we let $t \to 1$, and
complete the proof.
\ep
\section{Applications}
Previously we have examined how closed incompressible {\la} {\ms} interacts with {\mst}s of short curves
(Theorem ~\ref{main3} and Separation Lemma ~\ref{claim}). We now proceed to apply these techniques in the
settings of {\qfm}s and oriented closed {\htm}s that fiber over the circle, respectively.

\subsection{Multiplicity of {\ms}s in {\qfm}s}
When the {\cl} of the curve $\gamma \subset \Mscr$ satisfies the conditions in Theorem ~\ref{main1}, where $\Mscr$ is a {\qfm},
one would expect multiple {\ms}s around $\gamma$. This indeed the case, for instance, in Figure 3 below. This
is because there exist two closed incompressible surfaces $S_1$ and $S_2$ in $\Mscr\backslash\gamma$ which are not isotopic to each
other. Applying Theorem ~\ref{main3}, we produce two {\la} surfaces $T_1$ and $T_2$ that are not isotopic.

\begin{figure}[htbp]
\centering
\begin{minipage}[t]{1\linewidth}
\centering
\includegraphics[scale=0.75]{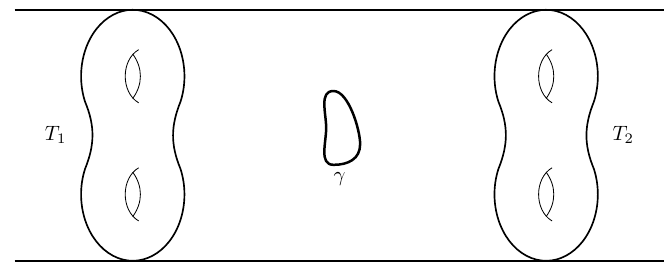}
\end{minipage}%
\caption{Two {\ms}s around a short curve.}
\end{figure}

We now make more precise this observation to the case of
multiple short (but unlinked) curves in a {\qfm}:
\bcor\label{coro2}
Let $\Gamma =\{\gamma_i\}_{i=1}^n$ be a collection of
mutually disjoint {\scg}s in a {\qfm} $\Mscr\cong{}S\times\R$,
each of {\cl} $\Lcal_i = \ell_i + \sqrt{-1}\,\theta_i$,
where $\ell_i>0$ and $\theta_i \in [-\pi,\pi)$, which satisfy:
\ben
\item $\ell_i < \varepsilon_{\Otal}(g)$;
\item $|\theta_i|/\sqrt{\ell_i} > \sqrt[4]{3\pi^2} \approx 2.33268$.
\een
If $\Sigma$ is an embedded closed incompressible surface
of genus $g\ge 2$ in $\Mscr\backslash\Gamma$
(which is homeomorphic to $S$), then there exists an embedded {\la} {\ms}
$T\subset\Mscr\backslash\Gamma$ isotopic to $\Sigma$.
Moreover, the {\qfm} $\Mscr$ contains at least $n+1$
distinct closed incompressible {\la} surfaces.
\ecor
\bp

By Theorem ~\ref{Mey}, and $\varepsilon_{\Otal}(g)<\varepsilon_2$,
the tubes $\mathbb{T}(\gamma_i)$ are mutually disjoint. Then the
first part of the corollary follows from the proof of Theorem ~\ref{main3}.

For the second part, with the real length condition
$\ell_i < \varepsilon_{\Otal}(g)$, for all $i = 1,2,\cdots, n$, the
collection $\Gamma$ is unlinked in the following sense (\cite{Ota03}):
there exists a homeomorphism between $\Mscr$ and $S\times\R$
such that each component of $\Gamma$ is
contained in one of the surfaces
$S\times \{i\}$, $ 1\le i\le n$. Now we count isotopy classes:
there are $n+1$ ways one can find closed incompressible surfaces
$\Sigma_{1},\ldots,\Sigma_{n+1}$ embedded
in $\Mscr\backslash\Gamma$ can separate
$\Gamma$ such that they are not isotopic to each other in
$\Mscr\backslash\Gamma$ (see Figure 4 for instance). For each
arrangement, we apply Theorem ~\ref{main3}, and then we find $n+1$
embedded closed incompressible {\la} surfaces $T_1,\ldots,T_{n+1}$
such that $T_{\alpha}$ is isotopic to $\Sigma_{\alpha}$ in
$\Mscr\backslash\Gamma$ for $\alpha=1,2,\ldots,n+1$. They are distinct
since they are not isotopic
pair-wisely.
\ep

\begin{figure}[htbp]
\centering
\begin{minipage}[t]{1\linewidth}
\centering
\includegraphics[scale=0.75]{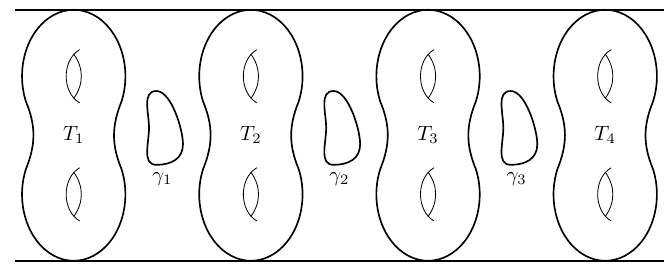}
\end{minipage}%
\caption{Minimal surfaces around multiple short curves.}
\end{figure}

\subsection{Proof of Theorem ~\ref{main1}}

We now move to our main interest: oriented closed {\htm}s that
fiber over the circle. Theorem \ref{main1} is proved by taking
algebraic limits in the {\qf} setting.
Now we will use Corollary ~\ref{coro2} to complete the proof of Theorem ~\ref{main1}, which we re-state here:

\begin{theorem1}
If an oriented closed {\htm} $\Mcal$ that fibers over the circle
with fiber $S$ contains a {\scg} whose {\cl}
$\Lcal = \ell + \sqrt{-1}\,\theta$ satisfies:
\ben
\item $\ell < \varepsilon_{\Otal}(g)$;
\item $|\theta|/\sqrt{\ell} > \sqrt[4]{3\pi^2} \approx 2.33268$,
\een
then for any positive integer $N$, there exists a {\qfm}
$\Mscr \cong S\times \R$ which contains
at least $N$ embedded closed incompressible {\la} {\ms}.
\end{theorem1}

\bp
Recall that $\Mcal$ is a closed {\htm} {\foc}, with
fiber $S$ closed surfaces of genus greater than one. We
consider a cyclic cover of $\Mcal$, ``unwrapping" the
circle direction. We denote this cover
$\Mcal_{\infty}\cong S\times\R$. Identifying $S$ with some
lift of the fiber, we obtain a discrete and faithful
representation $\rho: \pi_1(\Mcal_{\infty}) = \pi_1(S) \to \PSL$,
which is a {\ksg}.

Let $\gamma$ be a {\scg} on $\Mcal$ whose {\cl}
$\Lcal = \ell + \sqrt{-1}\,\theta$ satisfies
$0< \ell < \varepsilon_{\Otal}(g)$ and
$|\theta|/\sqrt{\ell}> \sqrt[4]{3\pi^2}$. By Proposition \ref{pro:unknotted},
the geodesic $\gamma$ can be lifted to a closed geodesic in $\Mcal_{\infty}$,
denoted by $\gamma$ too.
Let $\Phi$ be a deck transformation of $\Mcal_{\infty}$, then $\Mcal_{\infty}$
contains a sequence $\{\Phi^k(\gamma)\}_{k\in\mathbb{Z}}$ leaving every compact
subset in $\Mcal_{\infty}$ (see for
instance \cite{Min03}). This doubly degenerate {\htm} $\Mcal_{\infty}$
belongs to the Thurston boundary of
the deformation space, namely,
$\partial \QF(S) = \AH(S)\backslash\QF(S) = \overline{\QF_g(S)}\backslash\QF(S)$,
using \eqref{closure}. There is a sequence of {\qfg}s, each
representing a {\qfm} $\{\Mscr_i\}$, which
converges to the manifold $\Mcal_{\infty}$ algebraically
as $i \to \infty$.


\begin{figure}[htbp]
\centering
\begin{minipage}[t]{1\linewidth}
\centering
\includegraphics[scale=0.75]{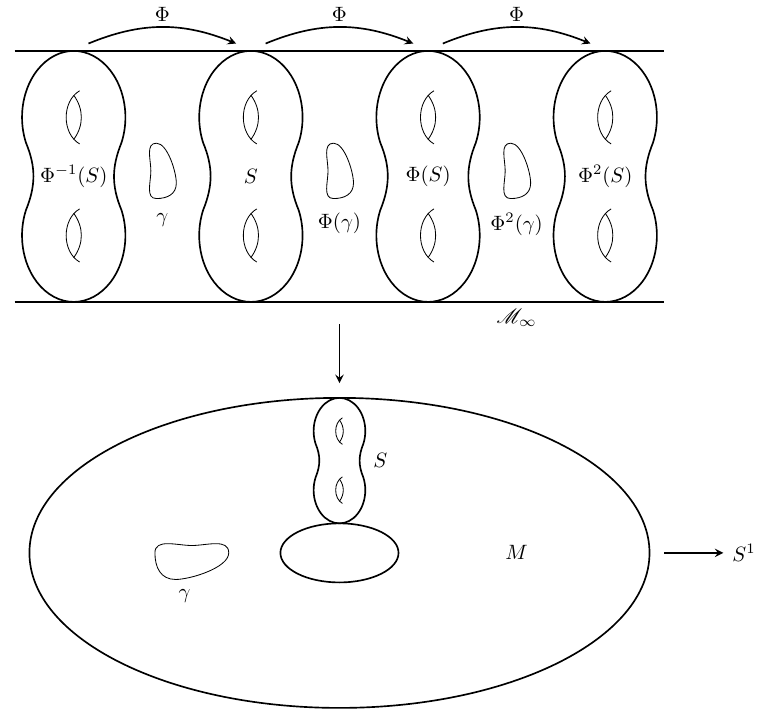}
\end{minipage}%
\caption{Cyclic cover for surface bundle {\foc}.}
\end{figure}

Since each element in a {\ksg} determines a geodesic, and
$\Mcal_{\infty}$, as a cyclic cover of $\Mcal$,
contains infinitely many hyperbolic geodesics
$\{\Phi^k(\gamma)\}_{k\in\mathbb{Z}}$, all having the
same {\cl} $\Lcal$. For any $N>0$, when $i$ is sufficiently
large, there is a {\qfm} $\Mscr_i$ in the
sequence such that it contains at least $N-1$ {\scg}s each
{\cl} satisfying the conditions in the statement.
We then apply Corollary ~\ref{coro2} to find that $\Mscr_i$
contains at least $N$ embedded closed incompressible {\la} surfaces.
\ep

\subsection{Proof of Theorem ~\ref{main2}}
In this subsection, we apply our estimates and cut-and-paste techniques developed in \S3 to prove the nonexistence
of a $C^2$-{\mf} for an oriented closed {\htm} that fibers over the circle, which we re-state here:

\begin{theorem2}
If an oriented closed {\htm} $\Mcal$ that fibers over the circle
with fiber $S$ contains a {\scg} whose {\cl}
$\Lcal = \ell + \sqrt{-1}\,\theta$ satisfies:
\ben
\item $\ell < \varepsilon_{\Otal}(g)$;
\item $|\theta|/\sqrt{\ell} > \sqrt[4]{3\pi^2} \approx 2.33268$,
\een
then $\Mcal$ does not admit a {\mf}.
\end{theorem2}

\bp
We proceed by contradiction. Suppose that the hyperbolic mapping
torus $\Mcal$ is foliated by {\ms}s all isotopic to a closed surface $S$.
We denote this $C^2$-foliation
by $\Fcal$. By theorems of Sullivan (\cite{Sul79}), Hass (\cite{Has86})
and Harvey-Lawson (\cite{HL82}),
all leaves of the foliation $\Fcal$ are of {\la} homologically,
and of the same area.

Since $\Mcal$ contains a {\scg} $\gamma$ whose {\cl} $\Lcal$ satisfies
conditions $\ell < \varepsilon_{\Otal}(g)$
and the inequality \eqref{ratio}, by Proposition \ref{pro:unknotted} the geodesic
$\gamma$ lies on an embedded surface isotopic to some fiber of $\Mcal$.
Our strategy is to prove:
(1) there exist leaves of $\Fcal$ which are disjoint from $\gamma$,
and (2) any such leaf of $\Fcal$ can't be very close to $\gamma$.

Firstly let's prove (1). By Proposition \ref{pro:unknotted},
there exists an embedded surface isotopic to the fiber which is
disjoint from $\gamma$, therefore we may apply Theorem \ref{main3}
to obtain a minimal surface $F'\subset\Mcal\setminus\gamma$ (with
respect to the {\hym} of $\Mcal$) which is isotopic to
the leaves of $\Fcal$.
Since each leaf of $\Fcal$ is a minimal surface,
$F'$ must be a leaf of $\Fcal$ by the {\maxp}.
This means there exist some leaves of $\Fcal$ which are disjoint from $\gamma$.

Secondly let's prove (2). There exist some other leaves of $\Fcal$ which must
intersect $\gamma$ since $\Fcal$ is a foliation of the manifold $\Mcal$,
therefore there exists a leaf $F\in\Fcal$ such that $F$ is disjoint from $\gamma$
and the (hyperbolic) distance between $F$ and $\gamma$ is less than $\delta$,
where $\delta$ is the constant defined by \eqref{delta}.
According to the argument in the proof of Theorem \ref{main3} and Lemma
\ref{claim}, this is impossible. Actually let $\Sigma$ be a component
of $F\cap\T$ which intersects the $\delta$-neighborhood of $\gamma$,
since  $F\cap\T=F\cap\mathbf{T}$,  where $\T$ is the {\mst}
(whose radius is still denoted by $r_0$) for this short geodesic
$\gamma \subset \Mcal$ and $\mathbf{T}=\T\setminus\gamma$. Then by Lemma \ref{lem:Lemma_5_2_HW15a},
$\Sigma$ is either a disk such that $\partial\Sigma$ is null-homotopic on $\partial\T$ or $\Sigma$ is an annulus
such that $\partial\Sigma$ consists of two isotopic essential slopes on $\partial\T$.

We claim that in both cases
\be\label{eq:large length II}
   \Length(\Sigma\cap\Nscr_{s}(\gamma))\ge{}2|\theta|\sinh(s)\ ,
   \quad \frac{r_{0}}{2}\le{}s\le{}r_0\ .
\ene
We consider two cases:
\ben
   \item $\Sigma$ is a disk. In this case, if there exists some
         $s'\in\left[\frac{r_{0}}{2},r_{0}\right]$ such that
         \eqref{eq:large length II} fails,
         then we have $\Length(\Sigma\cap\Nscr_{s'}(\gamma))<
         2|\theta|\sinh(s')\leq{}2\pi\sinh(s')$. We then apply the
         argument in the proof of Lemma~\ref{claim}, and find
         $\Sigma\cap\Nscr_{s'}(\gamma)$
         is disjoint from the $\delta$-neighborhood of $\gamma$,
         so is $\Sigma$. A contradiction.
   \item $\Sigma$ is an annulus. For any
         $s\in\left[\frac{r_{0}}{2},r_{0}\right]$, if
         $\Sigma\cap\Nscr_{s}(\gamma)$ only
         consists of disks, then \eqref{eq:large length II} is true
         according to the argument in the previous case. If
         $\Sigma\cap\Nscr_{s}(\gamma)$ contains at least one annulus,
         similar to the argument in the proof of Lemma ~\ref{claim},
         \eqref{eq:large length II} is also true.
\een
Therefore each component of $F\cap\T$ which intersects the $\delta$-neighborhood of
$\gamma$ must satisfies \eqref{eq:large length II}.

Then we have the area comparison as in \eqref{area-estimate}:
$$\Area(\Sigma)  > \Area(\partial\T).$$
This area estimate allows us to use cut-and-paste technique
again, namely, we can replace each component
of $F\cap\T$ that intersects the $\delta$-neighborhood of $\gamma$
by either an annulus or a disk on $\partial\T$,
to obtain a new surface $F''$ such that it
has less area than $F$, and isotopic to $F$. This is
impossible since $F$ is of the {\la}.

Now we have proved that if a leaf $F$ of the foliation
$\Fcal$ is disjoint from $\gamma$
(recall this kind of leaves always exist by the shrink-wrapping argument),
then it must be disjoint from the $\delta$-neighborhood of
$\gamma$, that is, we have proved (2).

But on the other hand, since $\Fcal$ is a foliation of $\Mcal$,
some leaves of $\Fcal$ must intersect $\gamma$. By the continuity of $\Fcal$,
some leaves of $\Fcal$ which are disjoint from $\gamma$ must be sufficiently
close to $\gamma$. This contradiction implies the non-existence of the minimal
foliation on $\Mcal$.

\ep
\subsection{Final remarks} We make several remarks.
\ben
\item
There are many other interesting work related the {\cl} with the geometry of {\htm}s, see for instance (\cite{Min99,Bre11,Mil14}).
\item
One may ask Question ~\ref{Th} for the existence of $C^0$ foliations. Our techniques rely on a Theorem of Sullivan (\cite{Sul79}) which requires
the taut foliation to be at least $C^2$.
\item
It is still unknown whether there exists any {\fhtm} which does admit a minimal foliation ($C^0$ or $C^2$). Theorem ~\ref{main2} provides sufficient
conditions for a negative answer to the existence of minimal foliations, and these conditions are verified in \S5.2 for many {\fhtm}s. One can
produce many more examples by varying the number of twists and the number of loops being twisted by Twister and SnapPy programs.
\item
One may further ask whether a {\htm} always admits a foliation of closed {\is}s of {\cmc}. We (\cite{HW13}) constructed a {\qfm} which does not
contains such a {\cmc} foliation, but the question remains open for many other cases of {\htm}s.
\een
\section{Appendix}
In this appendix, we want to explore the computational aspect of the ratio $\frac{|\theta|}{\sqrt{\ell}}$. Note that we will need the short curve $\gamma$
lying on an incompressible closed surface to work. Nevertheless we first give an upper bound for the ratio when $\ell$ is small, in terms of the genus
of the surface. In the second part, we use computer program Twister (\cite{Twister}) to produce some explicit examples. The first such examples were
produced for us by Saul Schleimer, to whom we are most grateful.

\subsection{An upper bound}
We present a proof of a statement by Thurston that if the geodesic is \emph{short}, then its rotation is \emph{small} in the following sense
(see the proof of Corollary 4.3 in \cite{Thu98}).

\begin{pro}[Thurston 1986]\label{ThurstonLemma}
Let $M\in\AH(S)$ be a complete {\htm}, here $S$ is an oriented
closed incompressible surface with genus $g(S)\geq{}2$. Let $\gamma\subset{}S$ be a {\scg} with {\cl} $\Lcal=\ell+\sqrt{-1}\,\theta$, where $\ell>0$ and
$-\pi\leq\theta<\pi$. If its real length is less than the Meyerhoff constant, namely, $\ell < \varepsilon_2 = \Ot(1) \approx 0.107071$, then
\be\label{smallrot}
   |\theta|<\frac{2\pi(g-1)}{\cosh(r_{0})-1}\ ,
\ene
where $r_{0}$ is the radius of the maximal solid torus of $\gamma$.
Furthermore, as $\ell\to{}0$, we have $\theta\to{}0$.
\end{pro}

\bp
According to \cite[Lemma 6.12]{MT98}, there exists a pleated surface
$f:\Sigma\to{}M$ whose pleating locus contains
$\gamma$, where $\Sigma\in\Tcal(S)$ is a hyperbolic surface. This pleated surface is nevertheless incompressible, so at least
one component of $\T\cap{}f(\Sigma)$ is an annulus whose core is the {\scg} $\gamma$.
The area of this annulus is greater than that of $\Ascr_a$ (see \eqref{annulus}) by the co-area formula, where
$a=|\theta|/\ell$.

The hyperbolic area of $\Sigma$ is $2\pi|\chi(\Sigma)|=4\pi(g-1)$, then we have
\beq
   \begin{aligned}
      4\pi(g-1)&>\Area(\T\cap{}f(\Sigma))\\
               &\geq\Area(\Ascr_{a})
                =2\int_{0}^{r_{0}}\sqrt{\ell^2\cosh^2(u)+
                 \theta^2\sinh^2(u)}\,du\\
               &\geq{}2|\theta|\int_{0}^{r_{0}}\sinh(u)du
                =2|\theta|(\cosh(r_{0})-1)\ ,
   \end{aligned}
\eeq
then we get \eqref{smallrot}.

Now we apply explicit formulas \eqref{tuberad} and \eqref{kappa} to examine the asymptotics. Since $r_{0}\to\infty$ as
$\ell\to{}0$, we have $\theta\to{}0$ as $\ell\to{}0$ by \eqref{smallrot}.
\ep
Note that if we assume $\ell < \varepsilon_{Otal}(g)$, then by Proposition ~\ref{pro:unknotted}, we may
replace the condition $\gamma \subset S$ in Proposition ~\ref{ThurstonLemma} to $\gamma \subset M$.

We examine the asymptotics in Proposition ~\ref{ThurstonLemma} more closely and find:
\bcor\label{rot}
Same assumption in above Proposition, we have, whenever $\ell$ is small enough,
\be\label{rot2}
   \frac{|\theta|}{\sqrt{\ell}} <
   4\pi\sqrt{\frac{4\pi}{\sqrt{3}}}\,(g-1).
\ene
\ecor
\bp
From \eqref{smallrot}, when $\ell < \varepsilon_2$, we have
\be\label{pre-exp}
\frac{|\theta|}{\sqrt{\ell}} < \frac{2\pi(g-1)}{\sqrt{\ell}(\cosh(r_{0})-1)}.
\ene
Given the explicit nature of $r_0$ in terms of $\sqrt{\ell}$ in \eqref{tuberad} and \eqref{kappa}, we expand the function
$\sqrt{\ell}(\cosh(r_{0})-1)$ in terms of $\sqrt{\ell}$ as follows:
\be\label{expan}
\sqrt{\ell}(\cosh(r_{0})-1) = \frac1b - \sqrt{\ell} - \frac{b}{24}\ell-\frac{353b^3}{5760}\ell^2 + o(\ell^{\frac52}),
\ene
where $b = \sqrt{\frac{4\pi}{\sqrt{3}}} \approx 2.69355$. Certainly from \eqref{expan}, we have, when $\ell >0$ is small enough,
\be\label{expan2}
\sqrt{\ell}(\cosh(r_{0})-1) > \frac{1}{2b},
\ene
Therefore we have
\be\label{ratio3}
\frac{|\theta|}{\sqrt{\ell}} < 4\pi b(g-1).
\ene
\ep
\br
For genus $g=2$, this upper bound is approximately $33.84815$, with a limit $16.92408$ as $\ell$ goes to zero. Clearly this
upper bound gets worse as $\ell$ goes from zero to $\varepsilon_2$. In comparison, in our main theorems, the lower bound for this ratio
we require is approximately $2.33268$. This propels us to look for {\fhtm}s with specific length spectra.
\er
\subsection{Examples via Twister and SnapPy}

Intuitively, in order to produce high rotational angle $\theta$, one may twist a loop in the {\tm} many times. In this subsection, we produce several explicit examples of {\fhtm}
which contains a closed curve with our conditions satisfied by computer programs.

In the first example, we let $S_2$ be a closed genus two surface, and $M$ be a mapping torus fibering over the circle. We consider $\gamma$ a short {\scg} on a fiber. Following
a suggestion of Saul Schleimer, we run Twister and SnapPy programs under the system Python. Note that for both main theorems, we look for the {\cl} to satisfy
$\ell < \Ot(g)$, and $\frac{|\theta|}{\sqrt{\ell}}  >  2.33268$. We explore the following:
\vskip 0.1in
\begin{verbatim}
>>> import twister
>>> S2 = twister.Surface('S_2')
>>> S2.info()
>>> A Twister surface of genus 2 with 0 boundary component(s)
Loops: a, b, c, d, e
Arcs:
>>> M = S2.bundle('b'*8 + 'cdea').high_precision()
\end{verbatim}
This produces a SnapPy manifold $M$ and we may ask the program to calculate its hyperbolic volume:
\vskip 0.1in
\begin{verbatim}
>>> M.volume()
7.991423345
\end{verbatim}
We now ask the program to specify the spectrum:
\vskip 0.1in
\begin{verbatim}
>>> M.length_spectrum(0.5)
\end{verbatim}
This returns a curve of {\cl} $0.1055786 + 0.84482566\sqrt{-1}$ in $M$, namely we have $\ell \approx 0.1055786$ and $\theta \approx 0.84482566$.
Though the ratio $\frac{|\theta|}{\sqrt{\ell}} \approx 2.60003$ is desirable, the real length is not short enough. We then choose to do the twists a few more times:
\vskip 0.1in
\begin{verbatim}
>>> M = S2.bundle('b'*25 + 'cdea').high_precision()
>>> M.volume()
8.142725
>>> M.length_spectrum(0.5)
\end{verbatim}
Now this returns a mapping torus $M$ of volume $8.142725$ and a closed curve of {\cl} $0.0098 + 0.25794\sqrt{-1}$. Now we have $\ell < \Ot(2) \approx 0.01515$, and the
ratio $\frac{|\theta|}{\sqrt{\ell}} \approx 2.60572$. This is an explicit example for both Theorems ~\ref{main1} and ~\ref{main2} for $g =2$.

One similarly can work with other surfaces and their bundles (over the circle) to obtain more examples of {\fhtm}s with our conditions satisfied. For instance, we find for the case of
genus 3:
\begin{verbatim}
>>> S3=twister.Surface((3,0))
>>> S3.info()
A Twister surface of genus 3 with 0 boundary component(s)
Loops: a0, b1, b2, b3, b4, b5, c
Arcs:
>>> M = S3.bundle('b1'*40 + 'a0b2b3b4b5c').high_precision()
>>> M.volume()
10.4355474
>>> M.length_spectrum(0.5)
\end{verbatim}
Now we obtain a short curve with $\ell = 0.00302$ and $\theta = 0.158958$, which yields the ratio $\frac{|\theta|}{\sqrt{\ell}} \approx 2.892537$. Note that when
$g=3$, Otal's constant $\Ot(3) \approx 0.00549389$. This is an example
for $g=3$ to satisfy both conditions in Theorems ~\ref{main1} and ~\ref{main2}.

We conclude with an example from the program Twister where fibering surfaces have genus 4 and contains a short curve with conditions in both Theorems satisfied. Note that
$\Ot(4) \approx 0.00280798$.
\begin{verbatim}
>>> S4 = twister.Surface((4,0))
>>> S4.info()
A Twister surface of genus 4 with 0 boundary component(s)
Loops: a0, b1, b2, b3, b4, b5, b6, b7, c
Arcs:
>>> M = S4.bundle('b1'*45 + 'a0b2b3b4b5b6b7c').high_precision()
>>> M.volume()
11.511256
>>> M.length_spectrum(0.5)
\end{verbatim}
One of the closed curves returned has {\cl} $0.002362+0.140781\sqrt{-1}$, which yields the ratio $\frac{|\theta|}{\sqrt{\ell}} \approx 2.8967$.

\bibliographystyle{amsalpha}
\bibliography{fiber}

\providecommand{\bysame}{\leavevmode\hbox to3em{\hrulefill}\thinspace}
\providecommand{\MR}{\relax\ifhmode\unskip\space\fi MR }
\providecommand{\MRhref}[2]{%
  \href{http://www.ams.org/mathscinet-getitem?mr=#1}{#2}
}
\providecommand{\href}[2]{#2}
\begin{thebibliography}{GMM01}

\bibitem[ACS06]{ACS06}
Ian Agol, Marc Culler, and Peter~B. Shalen, \emph{Dehn surgery, homology and
  hyperbolic volume}, Algebr. Geom. Topol. \textbf{6} (2006), 2297--2312.

\bibitem[Ago02]{Ago02}
Ian Agol, \emph{Volume change under drilling}, Geom. Topol. \textbf{6} (2002),
  905--916.

\bibitem[And83]{And83}
Michael~T. Anderson, \emph{Complete minimal hypersurfaces in hyperbolic
  {$n$}-manifolds}, Comment. Math. Helv. \textbf{58} (1983), no.~2, 264--290.

\bibitem[BB04]{BB04}
Jeffrey~F. Brock and Kenneth~W. Bromberg, \emph{On the density of geometrically
  finite {K}leinian groups}, Acta Math. \textbf{192} (2004), no.~1, 33--93.

\bibitem[BHS14]{Twister}
Mark Bell, Tracy Hall, and Saul Schleimer, \emph{Twister (computer software)},
  \url{https://bitbucket.org/Mark_Bell/twister/}, 2008--2014, Version 2.4.1.

\bibitem[Bon86]{Bon86}
Francis Bonahon, \emph{Bouts des vari\'et\'es hyperboliques de dimension
  {$3$}}, Ann. of Math. (2) \textbf{124} (1986), no.~1, 71--158.

\bibitem[Bre11]{Bre11}
William Breslin, \emph{Short geodesics in hyperbolic 3-manifolds}, Algebr.
  Geom. Topol. \textbf{11} (2011), no.~2, 735--745.

\bibitem[Bro07]{Bro07}
Kenneth Bromberg, \emph{Projective structures with degenerate holonomy and the
  {B}ers density conjecture}, Ann. of Math. (2) \textbf{166} (2007), no.~1,
  77--93.

\bibitem[CDGW]{SnapPy}
Marc Culler, Nathan~M. Dunfield, Matthias Goerner, and Jeffrey~R. Weeks,
  \emph{Snap{P}y, a computer program for studying the geometry and topology of
  $3$-manifolds}, Available at \url{http://snappy.computop.org} (2016).

\bibitem[CG06]{CG06}
Danny Calegari and David Gabai, \emph{Shrinkwrapping and the taming of
  hyperbolic 3-manifolds}, J. Amer. Math. Soc. \textbf{19} (2006), no.~2,
  385--446.

\bibitem[Chu68]{Chu68}
Vicki Chuckrow, \emph{On {S}chottky groups with applications to {K}leinian
  groups}, Ann. of Math. (2) \textbf{88} (1968), 47--61.

\bibitem[FHS83]{FHS83}
Michael Freedman, Joel Hass, and Peter Scott, \emph{Least area incompressible
  surfaces in {$3$}-manifolds}, Invent. Math. \textbf{71} (1983), no.~3,
  609--642.

\bibitem[GHW10]{GHW10}
Ren Guo, Zheng Huang, and Biao Wang, \emph{Quasi-{F}uchsian three-manifolds and
  metrics on {\TS}}, Asian J. Math. \textbf{14} (2010), no.~2, 243--256.

\bibitem[GMM01]{GMM01}
David Gabai, G.~Robert Meyerhoff, and Peter Milley, \emph{Volumes of tubes in
  hyperbolic 3-manifolds}, J. Differential Geom. \textbf{57} (2001), no.~1,
  23--46.

\bibitem[GW07]{GW07}
William Goldman and Richard Wentworth, \emph{Energy of twisted harmonic maps of
  {R}iemann surfaces}, In the tradition of {A}hlfors-{B}ers. {IV}, Contemp.
  Math., vol. 432, Amer. Math. Soc., Providence, RI, 2007, pp.~45--61.

\bibitem[Has86]{Has86}
Joel Hass, \emph{Minimal surfaces in foliated manifolds}, Comment. Math. Helv.
  \textbf{61} (1986), no.~1, 1--32.

\bibitem[Has05]{Has05}
\bysame, \emph{Minimal surfaces and the topology of three-manifolds}, Global
  theory of minimal surfaces, Clay Math. Proc., vol.~2, Amer. Math. Soc.,
  Providence, RI, 2005, pp.~705--724.

\bibitem[Has15]{Has15}
\bysame, \emph{Minimal fibrations and foliations of hyperbolic $3$-manifolds},
  Preprint (2015).

\bibitem[HL82]{HL82}
Reese Harvey and H.~Blaine Lawson, Jr., \emph{Calibrated foliations (foliations
  and mass-minimizing currents)}, Amer. J. Math. \textbf{104} (1982), no.~3,
  607--633.

\bibitem[HL12]{HL12}
Zheng Huang and Marcello Lucia, \emph{Minimal immersions of closed surfaces in
  hyperbolic three-manifolds}, Geom. Dedicata \textbf{158} (2012), 397--411.

\bibitem[HS88]{HS88}
Joel Hass and Peter Scott, \emph{The existence of least area surfaces in
  {$3$}-manifolds}, Trans. Amer. Math. Soc. \textbf{310} (1988), no.~1,
  87--114.

\bibitem[HW13]{HW13}
Zheng Huang and Biao Wang, \emph{On almost-{F}uchsian manifolds}, Trans. Amer.
  Math. Soc. \textbf{365} (2013), no.~9, 4679--4698.

\bibitem[HW15]{HW15b}
\bysame, \emph{Counting minimal surfaces in quasi-{F}uchsian manifolds}, Trans.
  Amer. Math. Soc. \textbf{367} (2015), 6063--6083.

\bibitem[HW17]{HW17}
\bysame, \emph{{Closed minimal surfaces in cusped hyperbolic three-manifolds}},
  Geom. Dedicata \textbf{189} (2017), no.~1, 17--37.

\bibitem[JK82]{JK82}
Troels J{\o}rgensen and Peter Klein, \emph{Algebraic convergence of finitely
  generated {K}leinian groups}, Quart. J. Math. Oxford Ser. (2) \textbf{33}
  (1982), no.~131, 325--332.

\bibitem[Mar74]{Mar74}
Albert Marden, \emph{The geometry of finitely generated {K}leinian groups},
  Ann. of Math. (2) \textbf{99} (1974), 383--462.

\bibitem[McM96]{Mcm96}
Curtis~T. McMullen, \emph{Renormalization and 3-manifolds which fiber over the
  circle}, Annals of Mathematics Studies, vol. 142, Princeton University Press,
  Princeton, NJ, 1996.

\bibitem[Mey87]{Mey87}
Robert Meyerhoff, \emph{A lower bound for the volume of hyperbolic
  {$3$}-manifolds}, Canad. J. Math. \textbf{39} (1987), no.~5, 1038--1056.

\bibitem[Mil14]{Mil14}
Christian Millichap, \emph{{Mutations and short geodesics in hyperbolic
  3-manifolds}}, ArXiv e-prints: 1406.6033 (2014).

\bibitem[Min99]{Min99}
Yair~N. Minsky, \emph{The classification of punctured-torus groups}, Ann. of
  Math. (2) \textbf{149} (1999), no.~2, 559--626.

\bibitem[Min03]{Min03}
\bysame, \emph{End invariants and the classification of hyperbolic
  3-manifolds}, Current developments in mathematics, 2002, Int. Press,
  Somerville, MA, 2003, pp.~181--217.

\bibitem[Mor82]{Mor82}
Hiroshi Mori, \emph{On surfaces of right helicoid type in {$\H^{3}$}}, Bol.
  Soc. Brasil. Mat. \textbf{13} (1982), no.~2, 57--62.

\bibitem[MT98]{MT98}
Katsuhiko Matsuzaki and Masahiko Taniguchi, \emph{Hyperbolic manifolds and
  {K}leinian groups}, Oxford Mathematical Monographs, The Oxford University
  Press, New York, 1998.

\bibitem[MY82a]{MY82}
William~H. Meeks, III and Shing~Tung Yau, \emph{The classical {P}lateau problem
  and the topology of three-dimensional manifolds. {T}he embedding of the
  solution given by {D}ouglas-{M}orrey and an analytic proof of {D}ehn's
  lemma}, Topology \textbf{21} (1982), no.~4, 409--442.

\bibitem[MY82b]{MY82b}
\bysame, \emph{The existence of embedded minimal surfaces and the problem of
  uniqueness}, Math. Z. \textbf{179} (1982), no.~2, 151--168.

\bibitem[Ota95]{Ota95}
Jean-Pierre Otal, \emph{Sur le nouage des g\'eod\'esiques dans les vari\'et\'es
  hyperboliques}, C. R. Acad. Sci. Paris S\'er. I Math. \textbf{320} (1995),
  no.~7, 847--852.

\bibitem[Ota96]{Ota96}
\bysame, \emph{Le th\'eor\`eme d'hyperbolisation pour les vari\'et\'es
  fibr\'ees de dimension 3}, Ast\'erisque (1996), no.~235.

\bibitem[Ota01]{Ota01}
\bysame, \emph{The hyperbolization theorem for fibered 3-manifolds}, SMF/AMS
  Texts and Monographs, vol.~7, American Mathematical Society, Providence, RI,
  2001, Translated from the 1996 French original by Leslie D. Kay.

\bibitem[Ota03]{Ota03}
\bysame, \emph{Les g\'eod\'esiques ferm\'ees d'une vari\'et\'e hyperbolique en
  tant que n\oe uds}, Kleinian groups and hyperbolic 3-manifolds ({W}arwick,
  2001), London Math. Soc. Lecture Note Ser., vol. 299, Cambridge Univ. Press,
  2003, pp.~95--104.

\bibitem[Rub07]{Rub07}
J.~H. Rubinstein, \emph{Problems around 3-manifolds}, Workshop on {H}eegaard
  {S}plittings, Geom. Topol. Monogr., vol.~12, Geom. Topol. Publ., Coventry,
  2007, pp.~285--298.

\bibitem[SU82]{SU82}
Jonathan Sacks and Karen~K. Uhlenbeck, \emph{Minimal immersions of closed
  {R}iemann surfaces}, Trans. Amer. Math. Soc. \textbf{271} (1982), no.~2,
  639--652.

\bibitem[Sul79]{Sul79}
Dennis Sullivan, \emph{A homological characterization of foliations consisting
  of minimal surfaces}, Comment. Math. Helv. \textbf{54} (1979), no.~2,
  218--223.

\bibitem[Sul85]{Sul85}
\bysame, \emph{Quasiconformal homeomorphisms and dynamics. {II}. {S}tructural
  stability implies hyperbolicity for {K}leinian groups}, Acta Math.
  \textbf{155} (1985), no.~3-4, 243--260.

\bibitem[SY79]{SY79}
Richard Schoen and Shing~Tung Yau, \emph{Existence of incompressible minimal
  surfaces and the topology of three-dimensional manifolds with nonnegative
  scalar curvature}, Ann. of Math. (2) \textbf{110} (1979), no.~1, 127--142.

\bibitem[Thu80]{Thu80}
William~P. Thurston, \emph{The geometry and topology of three-manifolds},
  Princeton University, 1980.

\bibitem[Thu82]{Thu82}
\bysame, \emph{Three-dimensional manifolds, {K}leinian groups and hyperbolic
  geometry}, Bull. Amer. Math. Soc. (N.S.) \textbf{6} (1982), no.~3, 357--381.

\bibitem[Thu98]{Thu98}
\bysame, \emph{Hyperbolic structures on 3-manifolds, {II}: Surface groups and
  3-manifolds which fiber over the circle}, ArXiv Mathematics e-prints (1998).

\bibitem[Tuz93]{Tuz93}
A.~A. Tuzhilin, \emph{Global properties of minimal surfaces in {${\bf R}^3$}
  and {${\bf H}^3$} and their {M}orse type indices}, Minimal surfaces, Adv.
  Soviet Math., vol.~15, Amer. Math. Soc., Providence, RI, 1993, pp.~193--233.

\bibitem[Uhl83]{Uhl83}
Karen~K. Uhlenbeck, \emph{Closed minimal surfaces in hyperbolic
  {$3$}-manifolds}, Seminar on minimal submanifolds, Ann. of Math. Stud., vol.
  103, Princeton Univ. Press, Princeton, NJ, 1983, pp.~147--168.

\bibitem[Wan12]{Wan12}
Biao Wang, \emph{Minimal surfaces in quasi-{F}uchsian 3-manifolds}, Math. Ann.
  \textbf{354} (2012), no.~3, 955--966.

\bibitem[Wan16]{Wan16}
\bysame, \emph{Stability of catenoids and helicoids in hyperbolic space}, Asian
  J. Math. \textbf{to appear} (2016).

\bibitem[Wie77]{Wie77}
Norbert~J. Wielenberg, \emph{Discrete {M}oebius groups: fundamental polyhedra
  and convergence}, Amer. J. Math. \textbf{99} (1977), no.~4, 861--877.

\end{thebibliography}
\end{document}